\documentclass[titlepage,11pt]{article}
\usepackage{latexsym}
\usepackage{amsfonts}
\usepackage{amsmath}
\usepackage{mathtools}
\usepackage{tikz}
\usepackage{float}
\usepackage{comment}
\usetikzlibrary{arrows}
\usetikzlibrary{decorations.markings}
\usetikzlibrary{decorations.pathreplacing,calligraphy}
\usetikzlibrary{positioning,arrows.meta}

\newcommand\blackslug{\hbox{\hskip 1pt \vrule width 4pt height 8pt depth 1.5pt
		\hskip 1pt}}
\newcommand\bbox{\hfill \quad \blackslug \bigbreak}
\def\DD{\hbox{-}}
\def\CC{\hbox{-}\cdots\hbox{-}}
\def\LL{,\ldots,}
\newcommand{\vare}{\varepsilon}

\newcommand{\cupcup}{\cup \cdots\cup}

\def\dist{\operatorname{dist}}
\newcommand{\mac}{\mathcal}

\DeclarePairedDelimiter\abs{\lvert}{\rvert}%

\title{Asymptotic structure. V. The coarse Menger conjecture in bounded path-width\\ 
}
\author{
	Alex Divoux\\
	Princeton University,\\ Princeton, NJ 08544, USA
	\and
	Tung Nguyen\thanks{Part of this work was conducted while the second author was at Princeton University and was supported by a Porter Ogden Jacobus Fellowship. Currently supported by a Titchmarsh Research Fellowship and a Christ Chruch Research Centre Grant.}\\
	University of Oxford,\\ Oxford, UK
	\and
	Alex Scott\thanks{Supported by EPSRC grant EP/X013642/1}\\
	University of Oxford, \\
	Oxford, UK
	\and
	Paul Seymour\thanks{Supported by AFOSR grant
		FA9550-22-1-0234, and by NSF grant DMS-2154169.}\\
	Princeton University,\\ Princeton, NJ 08544, USA}

\date{}

\newtheorem{thm}{}[section]

\newcommand{\Proof}{\noindent{\bf Proof.}\ \ }

\begin{document}
	\maketitle
	\begin{abstract}
		Menger's theorem tells us that if $S,T$ are sets of vertices in a graph $G$, then (for $k\ge0$) either there are $k+1$ vertex-disjoint paths between $S$ and $T$, or there is a set of $k$ vertices separating $S$ and $T$.  But what if we want the paths to be far apart, say at distance at least $c$?  
		One might hope that we can find either $k+1$ paths pairwise far apart, or $k$ sets of bounded radius that separate $S$ and $T$, where 
		the bound on the radius is some $\ell$ that depends only on $k,c$
		(the ``coarse Menger conjecture'').
		The last three authors showed in an earlier paper that this is false for all $k\ge 2$ and $c\ge3$, by constructing a sequence of finite graphs giving counterexamples
		for larger and larger values of $\ell$ with $k=2$ and $c=3$. These counterexamples contained subdivisions of uniform binary trees with 
		arbitrarily large depth as subgraphs, and so had unbounded path-width.

		Here we show that, if $H$ is a graph that can be drawn in the plane such that each region shares a vertex with the infinite region,
		then the coarse Menger conjecture is true for all graphs not containing $H$ as a minor. Consequently, the conjecture is true for 
		all graphs with bounded path-width (by taking $H$ to be a sufficiently large tree), and it is true for series-parallel graphs 
		(by taking $H=K_4$).  The first is somewhat surprising, since the conjecture is false for bounded tree-width.
		
	\end{abstract}
	
	\section{Introduction}
	
	Let $S,T$ be sets of vertices of a graph $G$. (In this paper, all graphs are finite and have no loops or multiple edges.)
	Menger's theorem~\cite{menger} tells us that either there are $k+1$ pairwise vertex-disjoint paths between $S$ and $T$, or there is a 
	set $X$ of at most $k$ vertices such that every $S$-$T$ path in $G$ meets $X$. But what if we want the paths to be pairwise far apart? In this case,  the question is much harder.  Bienstock~\cite{bienstock} showed that it is NP-hard to decide whether, 
	given four vertices $s_1,s_2,t_1,t_2$ of a graph $G$, there are two paths between 
	between $\{s_1,s_2\}$ and $\{t_1,t_2\}$ that have distance $\ge 2$, that is, they are vertex-disjoint and there is no edge joining them. 
	This was recently extended by Balig\'acs and MacManus~\cite{bm}, who showed the same thing for distance $\ge c$, for each $c\ge 3$.

	Since the problem is NP-complete, one would not expect to find a necessary and sufficient condition for the existence of $k+1$ $S$-$T$ 
	paths at distance at least $c$; but still one could hope for some sort of obstruction that is necessary for excluding $k+1$ $S$-$T$ paths at 
	distance at least $c$, and sufficient for excluding $k+1$ $S$-$T$ paths at distance more than some larger number depending on $k,c$.
	Two groups of researchers, Albrechtsen, Huynh, Jacobs, Knappe and Wollan~\cite{wollan}, and independently Georgakopoulos and Papasoglu~\cite{agelos}, proposed such a statement:
	\begin{thm}\label{conj}
		{\bf Coarse Menger Conjecture:} For all integers $k\ge 0$ and $c\ge1$ there exists $\ell\ge 0$ with the following property.
		Let $G$ be a graph and let $S,T\subseteq V(G)$.  Then either
		\begin{itemize}
			\item there are $k+1$ paths between $S,T$, pairwise at distance at least $c$; or
			\item there is a set $X\subseteq V(G)$
			with $|X|\le k$ such that every path between $S,T$ contains a vertex at distance at most $\ell$ from some member of $X$.
		\end{itemize}
	\end{thm}
	Both groups showed that this is true for $k=1$, 
	and 
	Gartland, Korhonen and Lokshtanov~\cite{gkl} and Hendrey, Norin, Steiner, and Turcotte~\cite{hnst} proved it for bounded degree graphs 
	when $c = 2$.
	However, 
	the last three authors showed in~\cite{counterex} that the coarse Menger conjecture is false for all $k\ge 2$, for any fixed $c\ge 3$. Indeed, it remains false even if we weaken the bound $|X|\le k$
	in the second bullet to $|X|\le m$, where $m$ is any constant depending on $k,c$~\cite{AS4}.
	
	Thus, we need to lower our sights a little, and one way to do so is to work in restricted classes of graphs. The counterexamples
	of~\cite{counterex} have unbounded genus, and have unbounded 
	``path-width'' (defined in the next section).  
	It might be true that the coarse Menger conjecture holds for graphs of  bounded genus, but this is open; see~\cite{AS6} for some progress in this direction. Our result implies that the coarse Menger conjecture is true
	for graphs of bounded path-width. 
	
	More exactly, it implies:
	\begin{thm}\label{mainthm}
		Let $k,d\ge 0$ and $c\ge1$ be integers. Then there exists $\ell\ge 0$, such that for every graph $G$ with path-width at most $d$, and all $S,T\subseteq V(G)$, either:
		\begin{itemize}
			\item there are $k+1$ paths between $S,T$, pairwise at distance at least $c$; or
			\item there is a set $X\subseteq V(G)$
			with $|X|\le k$ such that every path between $S,T$ contains a vertex at distance at most $\ell$ from some member of $X$.
		\end{itemize}
	\end{thm}
	(For fixed $k,d$, the number $\ell$ depends linearly on $c$, as we will see.)
	We will deduce in the conclusion that the coarse Menger conjecture 
	also holds for graphs of ``bounded coarse path-width''.
	Curiously, the coarse Menger conjecture is {\em not} true for graphs of bounded tree-width, since the counterexamples
	of~\cite{counterex} have tree-width six.
	
	Our main result is more general than this. 
	Let us say a graph $H$ is {\em near-outerplanar} if it is planar, and has a planar drawing such that every region is incident with a vertex of the infinite region (see Figure \ref{fig:nearouterplanar}). Thus, outerplanar graphs are near-outerplanar, but so are many other graphs, such as wheels.
	
	\begin{figure}[H]
		\centering
		
		\begin{tikzpicture}[auto=left]
			\tikzstyle{every node}=[inner sep=1.5pt, fill=black,circle,draw]
			\def\r{2}
			\node (0) at ({\r*cos(0)},{\r*sin(0)}) {};
			\node (1) at ({\r*cos(30)},{\r*sin(30)}) {};
			\node (2) at ({\r*cos(60)},{\r*sin(60)}) {};
			\node (3) at ({\r*cos(90)},{\r*sin(90)}) {};
			\node (4) at ({\r*cos(120)},{\r*sin(120)}) {};
			\node (5) at ({\r*cos(150)},{\r*sin(150)}) {};
			\node (6) at ({\r*cos(180)},{\r*sin(180)}) {};
			\node (7) at ({\r*cos(210)},{\r*sin(210)}) {};
			\node (8) at ({\r*cos(240)},{\r*sin(240)}) {};
			\node (9) at ({\r*cos(270)},{\r*sin(270)}) {};
			\node (10) at ({\r*cos(300)},{\r*sin(300)}) {};
			\node (11) at ({\r*cos(330)},{\r*sin(330)}) {};
			
			\node (12) at (0,0) {};
			\node (13) at (1,1) {};
			\node (14) at (-.7,0) {};
			\node (15) at (-.3,-1.3) {};
			\node (16) at (.7,-.7) {};
			\node (17) at (-.8,.7) {};
			\node (18) at (-.3,1.3) {};
			\node (19) at (-1.5,.4) {};
			\node (20) at (1.2,-.3) {};
			\node (21) at (1.2,.3) {};
			
			\foreach \from/\to in {0/1,1/2,2/3,3/4,4/5,5/6,6/7,7/8,8/9,9/10,10/11,11/0, 12/13,12/14,12/15,12/16,13/1,13/2,14/7,14/17,
				15/8,15/9,16/20,16/10,17/18,17/19,18/3,18/4,19/5,19/6,20/0,20/11,17/4,5/17,12/21,21/1,21/0}
			\draw [-] (\from) -- (\to);
			
		\end{tikzpicture}
		
		\caption{A near-outerplanar graph.}
		\label{fig:nearouterplanar}
	\end{figure}
	
	We will prove that:
	\begin{thm}\label{outerplanethm}
		For every near-outerplanar graph $H$, and all integers $k\ge 0$ and $c\ge1$, there exists $\ell\ge 0$, such that for every graph $G$ 
		not containing $H$ as a minor, and all $S,T\subseteq V(G)$, either:
		\begin{itemize}
			\item there are $k+1$ paths between $S,T$, pairwise at distance at least $c$; or
			\item there is a set $X\subseteq V(G)$
			with $|X|\le k$ such that every path between $S,T$ contains a vertex at distance at most $\ell$ from some member of $X$.
		\end{itemize}
	\end{thm}
	(We could assume that $k\ge 2$ if we wanted, because the coarse Menger conjecture is known to be true for
	$k=1$~\cite{wollan, agelos,counterex}, but there is no need.)
	In particular, every tree $H$ is near-outerplanar, and so the coarse Menger conjecture is true for graphs that do not contain $H$ as a minor; and this implies \ref{mainthm}, as we discuss in the next section. Moreover, the graph $K_4$ is near-outerplanar, and so 
	the coarse Menger conjecture is true for series-parallel graphs, since they are the graphs with no $K_4$ minor. (This was an open question.)

	\section{Subdivisions and path-width}
	
	In this paper, for $d\ge 1$, the ``uniform binary tree of depth $d+1$'' is the tree $H$ such that for some $r\in V(H)$ (the ``root''), $r$ has degree two, all other vertices
	have degree one or three, and every vertex of degree one has distance exactly $d$ from $r$. Thus, $H$ has $2^{d+1}-1$ vertices. 
	We denote this tree by $H_{d}$. 
	
	If $H$ is a graph, a {\em subdivision} of $H$ is a graph obtained from $H$ by replacing each of its edges by a path of length at
	least one joining the same pair of vertices, where these paths are pairwise vertex-disjoint except for their ends. For $\ell\ge 1$,
	let us say an
	{\em $(<\ell)$-subdivision} of $H$ is a subdivision obtained by replacing each edge by a path of length $\le \ell$ (and at least one).
	(This is to be consistent with
	the standard term ``1-subdivision'', which means replacing each edge with a path
	of length two.)
	
	Let us define path-width. A graph $G$ has {\em path-width} at most $d$ if and only if there is a sequence $W_1\LL W_n$ of subsets of its vertex set, satisfying:
	\begin{itemize}
		\item $|W_i|\le d+1$ for $1\le i\le n$;
		\item $G[W_1]\cupcup G[W_n] = G$; and
		\item $W_i\cap W_k\subseteq W_j$ for $1\le i\le j\le k\le n$.
	\end{itemize}

	We do not really need this definition. The only thing about bounded path-width that concerns us is a theorem of Robertson and Seymour~\cite{GM1}:
	\begin{thm}\label{GM1}
		For every integer $\delta\ge 1$, there exists $k$, such that every graph that contains no subdivision of
		$H_\delta$ as a subgraph has path-width at most $k$; and conversely, every graph that contains a
		subdivision of $H_\delta$ as a subgraph has path-width at least $\delta/2$.
	\end{thm}
	Thus, knowing that there is a bound on path-width is the same as knowing that for some $\delta$, no subgraph is a  subdivision of 
	$H_\delta$.
	Indeed, in this paper it is more natural to work with the ``excluded tree subdivision'' version directly, rather than working with path-width. 
	So we could reformulate \ref{mainthm} as follows:
	\begin{thm}\label{treeform}
		For all integers $k,\delta,c\ge 0$ there exist $\ell\ge 0$, with the following property. Let $G$ be a graph that contains no 
		subdivision of $H_\delta$ as a subgraph, and let $S,T\subseteq V(G)$. Then either
		\begin{itemize}
			\item there are $k+1$ paths between $S,T$, pairwise at distance greater than $c$; or
			\item there is a set $X\subseteq V(G)$
			with $|X|\le k$ such that every path between $S,T$ contains a vertex at distance at most $\ell$ from some member of $X$.
		\end{itemize}
	\end{thm}
	Indeed, one can easily modify the proof we will give (although we will not do so), to show that the same conclusion holds even if we 
	just assume that $G$ contains no $(<\ell)$-subdivision of $H_\delta$.

	\section{Near-outerplanar graphs and the bordered binary tree.}
	
	But, as we said, our main result is more general than \ref{treeform}.
	For $\delta\ge 1$, the graph $H_\delta$ is a uniform binary tree of depth $\delta+1$.
	Let us add a path running though the leaves in this order, as in Figure \ref{fig:border}. We call the graph that results $H^+_\delta$. 
	The vertex of degree two at the top of the figure is its {\em root}. 
	
	\begin{figure}[h!]
		\centering
		
		\begin{tikzpicture}[scale=1/2,auto=left]
			
			\tikzstyle{every node}=[inner sep=1.5pt, fill=black,circle,draw]
			\node (v2) at (2,0) {};
			\node (v3) at (3,0) {};
			\node (v4) at (4,0) {};
			\node (v5) at (5,0) {};
			\node (v6) at (6,0) {};
			\node (v7) at (7,0) {};
			\node (v8) at (8,0) {};
			\node (v9) at (9,0) {};
			\node (v10) at (10,0) {};
			\node (v11) at (11,0) {};
			\node (v12) at (12,0) {};
			\node (v13) at (13,0) {};
			\node (v14) at (14,0) {};
			\node (v15) at (15,0) {};
			\node (v16) at (16,0) {};
			\node (v17) at (17,0) {};
			\node (v18) at (18,0) {};
			\node (v19) at (19,0) {};
			\node (v20) at (20,0) {};
			\node (v21) at (21,0) {};
			\node (v22) at (22,0) {};
			\node (v23) at (23,0) {};
			\node (v24) at (24,0) {};
			\node (v25) at (25,0) {};
			\node (v26) at (26,0) {};
			\node (v27) at (27,0) {};
			\node (v28) at (28,0) {};
			\node (v29) at (29,0) {};
			\node (v30) at (30,0) {};
			\node (v31) at (31,0) {};
			\node (v32) at (32,0) {};
			\node (v33) at (33,0) {};

			\draw[-] (v2)--(v33);
			
			\node (u2) at (2.5,1) {};
			\node (u4) at (4.5,1) {};
			\node (u6) at (6.5,1) {};
			\node (u8) at (8.5,1) {};
			\node (u10) at (10.5,1) {};
			\node (u12) at (12.5,1) {};
			\node (u14) at (14.5,1) {};
			\node (u16) at (16.5,1) {};
			\node (u18) at (18.5,1) {};
			\node (u20) at (20.5,1) {};
			\node (u22) at (22.5,1) {};
			\node (u24) at (24.5,1) {};
			\node (u26) at (26.5,1) {};
			\node (u28) at (28.5,1) {};
			\node (u30) at (30.5,1) {};
			\node (u32) at (32.5,1) {};
			
			\node (t3) at (3.5,2) {};
			\node (t7) at (7.5,2) {};
			\node (t11) at (11.5,2) {};
			\node (t15) at (15.5,2) {};
			\node (t19) at (19.5,2) {};
			\node (t23) at (23.5,2) {};
			\node (t27) at (27.5,2) {};
			\node (t31) at (31.5,2) {};
			
			\draw (u2) -- (t3)--(u4);
			\draw (u6) -- (t7)--(u8);
			\draw (u10) -- (t11)--(u12);
			\draw (u14) -- (t15)--(u16);
			\draw (u18) -- (t19)--(u20);
			\draw (u22) -- (t23)--(u24);
			\draw (u26) -- (t27)--(u28);
			\draw (u30) -- (t31)--(u32);
			
			\node (s5) at (5.5,3) {};
			\node (s13) at (13.5,3) {};
			\node (s21) at (21.5,3) {};
			\node (s29) at (29.5,3) {};
			
			\draw (t3) -- (s5)--(t7);
			\draw (t11) -- (s13)--(t15);
			\draw (t19) -- (s21)--(t23);
			\draw (t27) -- (s29)--(t31);
			
			\node (r9) at (9.5,4) {};
			\node (r25) at (25.5,4) {};
			
			\draw (s5) -- (r9)--(s13);
			\draw (s21) -- (r25)--(s29);
			\node (q17) at (17.5,5) {};
			\draw (r9) -- (q17)--(r25);
			
			\foreach \from/\to in {v2/u2,u2/v3,v4/u4,u4/v5,v6/u6,u6/v7, v8/u8,u8/v9,v10/u10,u10/v11,v12/u12,u12/v13,v14/u14,u14/v15,v16/u16,
				u16/v17,v18/u18,u18/v19,v20/u20,u20/v21,v22/u22,u22/v23,v24/u24,u24/v25,v26/u26,u26/v27,v28/u28,u28/v29,v30/u30,u30/v31,v32/u32,u32/v33}
			\draw [-] (\from) -- (\to);
			
		\end{tikzpicture}
		
		\caption{The graph $H_5^+$.} \label{fig:border}
	\end{figure}
	The added path in the graph of the figure runs through all 32 leaves of $H_5$; let us number them $v_1\LL v_{32}$ in order. There are 
	many other orders that would give an isomorphic graph: for instance a path running though the leaves in the order 
	$v_2,v_1,v_3,v_4\LL v_{32}$ would give a graph isomorphic to $H_5^+$, as would the order $v_3,v_4,v_1,v_2,v_5\LL v_{32}$. But not every order works; let us say a linear order of the leaves of $H_\delta$ is {\em normal} if adding a path through the leaves in that order
	gives a copy of $H^+_\delta$. 
	
	Evidently, $H_\delta^+$ is near-outerplanar, and it is an easy exercise to show that every near-outerplanar graph is isomorphic to a 
	minor of $H_\delta^+$ for sufficiently large $\delta$. 
	To see this, let $H$ be near-outerplanar, drawn such that every region is incident with a vertex on the infinite region.
	\begin{itemize}
		\item By adding an edge between two components if possible, and repeating, we may assume that $H$ is connected.
		\item By adding an edge between $ac$ if $a\DD b\DD c$ is a path of the infinite boundary and $b$ is a cut-vertex separating $a,c$,
		and repeating, we may assume the boundary of the infinite region is a cycle $C$.
		\item If some vertex of $C$ has degree at least four, we may uncontract an edge (making $C$ longer), ``splitting'' this vertex into 
		two adjacent vertices of smaller degree. By repeating, we may assume that every vertex of $C$ has degree at most three.
		\item Hence, $H\setminus E(C)$ is a forest, and by subdividing and adding edges we may assume it is a tree $T$, and all leaves of $T$ are in $C$.
		\item We may assume this tree has maximum degree at most three, by splitting vertices into two adjacent vertices of smaller degree as before.
		\item And then the result is clear.
	\end{itemize}
	We will show that:
	\begin{thm}\label{border2} 
		For all integers $\delta\ge 5$, $k\ge 0$,
		and $c\ge 0$, let $G$ be a graph with no $H_\delta^+$ minor,
		and let $S,T\subseteq V(G)$. Then either
		\begin{itemize}
			\item there are $k+1$ paths between $S,T$, pairwise at distance greater than $c$; or
			\item there is a set $X\subseteq V(G)$
			with $|X|\le k$ such that every path between $S,T$ contains a vertex at distance at most $(k\delta)^{4k\delta} c$ from some member of $X$.
		\end{itemize}
	\end{thm}
	This strengthens \ref{treeform}, and will be our main result.
	We remark that here we are asking for paths at distance $>c$ rather than $\ge c$ as in \ref{mainthm}; we find this form slightly
	more convenient. We are assuming $\delta\ge 5$ for numerical reasons.

	The main purpose of this section is to prove that if we have a sufficiently large uniform binary tree with a path running through many of its leaves 
	in any order, then $G$ contains $H_\delta^+$ as a minor. 
	
	\begin{thm}\label{gathering}
		Let $\delta,h\ge 1$; let $T$ be a subgraph of $G$ isomorphic to a subdivision of $H_{h}$, with root $r$, and let $P$ be a path of $G$ that contains 
		more than $(2h)^{\delta-1}$ leaves of $T$,
		and contains no other vertices of $T$. Then there is a subtree $J$ of $T$, such that:
		\begin{itemize}
			\item $J$, rooted at $s$, is isomorphic to a subdivision of the 
			rooted tree $H_{\delta}$, where $s$ is the vertex of $J$ closest in $T$ to $r$;
			\item all leaves of $J$ are leaves of $T$; and
			\item $P$ contains all the leaves of $J$ in normal order, and contains no other vertices of $J$. 
		\end{itemize}
		Consequently $G$ contains a subdivision of $H_{\delta}^+$ as a subgraph.
	\end{thm}
	\Proof We may assume that $T$ is isomorphic to $H_{h}$, and we proceed by induction on $h$. If $\delta=1$, then $P$ contains at least two leaves of $H_h$ and the claim holds, so we assume 
	that $\delta\ge 2$. If $h=1$, then $T$ has only two leaves, and yet by hypothesis $P$ contains more than $2^{\delta-1}$ of them, which
	is impossible. Thus we may assume that $h\ge 2$ and the result holds for $h-1$.
	Let $r_1,r_2$ be the neighbours of $r$ in $T$, and let $T_1,T_2$ be the components of $T\setminus r$, where $r_i\in V(T_i)$
	for $i = 1,2$. Suppose first that $|V(P\cap T_2)|\le 2(2h)^{\delta-2}$. Then 
	$$|V(P\cap V(T_1)|> (2h)^{\delta-1}-2(2h)^{\delta-2}=(2h-2)(2h)^{\delta-2}\ge (2h-2)^{\delta-1},$$
	and the result follows from the inductive hypothesis applied to $T_1$ and $P$. So we may assume that $|V(P\cap T_2)|> 2(2h)^{\delta-2}$,
	and similarly $|V(P\cap T_1)|> 2(2h)^{\delta-2}$. 
	
	Choose an edge $e$ of $P$ such that both components $P_1,P_2$ of $P\setminus e$
	contain more than $(2h)^{\delta-2}$ vertices of $T_1$. Then one of $P_1,P_2$ contains more than $(2h)^{\delta-2}$ vertices of $T_2$,
	(since $|V(P\cap T_2)|> 2(2h)^{\delta-2}$)
	and we may assume that it is $P_2$. From the inductive hypothesis, for $i = 1,2$, there is a subtree $J_i$ of $T_i$, satisfying
	the three bullets above, with $T,r,\delta,J$ replaced by $T_i, r_i, \delta-1, J_i$. Let $J$ be the union of $J_1,J_2$ and the path of $T$ 
	between $r_1,r_2$ (which passes through $r$). Then $J$, rooted at $r$, is isomorphic to a subdivision of $H_{\delta}$. Moreover, all leaves 
	of $J$ are leaves of one of $J_1,J_2$ and so are leaves of $T$; and $P$ contains all the leaves of $J$, and no 
	other vertices of $J$. But all vertices of $J$ in $J_i$ belong to $P_i$ for $i = 1,2$, and since $P$ passes through the leaves 
	of $J_1$ and of $J_2$
	in normal order, it also passes through the leaves of $J$ in normal order. This proves \ref{gathering}.~\bbox

	\section{A key lemma}
	
	If we contract an edge of a graph, then distances do not change by much, but if we delete an edge or a vertex, they might change 
	considerably.
	In this section, we prove a lemma that allows us to bypass this problem to some extent,  in graphs excluding subdivisions of some 
	$H_\delta^+$.
	The proofs in this section are 
	the only places in the paper where we use the hypothesis about subdivisions of $H_\delta^+$.
	
	If $X$ is a vertex of a graph $G$, or
	a subset of the vertex set of $G$, or a subgraph of $G$, and the same for $Y$, then $\dist_G(X,Y)$ denotes the
	distance in $G$
	between $X,Y$, that is, the number of edges in the shortest path of $G$ with one end in $X$ and the other in $Y$. (If no path exists we set $\dist_G(X,Y) = \infty$.)

	For convenience, when $Y\subseteq V(G)$, we sometimes write $\dist_Y(u,v)$ for $\dist_{G[Y]}(u,v)$. 
	If $G$ is a graph, a {\em tie-breaker weighting} is a map $\lambda:E(G)\to \mathbb R$ with the following properties, where $\lambda(X)$ means 
	$\sum_{e\in X}\lambda(e)$:
	\begin{itemize}
		\item if $X,Y\subseteq E(G)$ and $\lambda(X)\le \lambda(Y)$ then $|X|\le |Y|$;
		\item if $X,Y\subseteq E(G)$ are distinct, then $\lambda(X)\ne \lambda(Y)$. 
	\end{itemize}
	To obtain such a map, one could define $\lambda(e)=1+\vare(e)$, where the numbers $\vare(e)$ are independent and very small.
	If $\lambda$ is a tie-breaker weighting, we say a path $P$ is a {\em $\lambda$-geodesic} if there is no path $Q$ with the same ends as $P$ and with
	$\lambda(E(Q))<\lambda(E(P))$. For all $u,v$ in the same component, there is a unique $\lambda$-geodesic between $u,v$. 
	
	Let $Z\subseteq V(G)$, and let $H\subseteq G$ be a tree with $|V(H)|\ge 2$, with all its leaves in $Z$ and no other vertices in $Z$.
	We say that $H$ is {\em $Z$-leaved}. 
	
	\begin{thm}\label{newkey}
		Let $\delta,\ell\ge 1$, let $G$ be a graph, and let $Z\subseteq V(G)$  such that no subgraph is a $Z$-leaved $(<2^{\delta-1}\ell)$-subdivision of $H_{\delta}$. 
		Then there exists $Y\supseteq Z$ with the following properties:
		\begin{itemize}
			\item every vertex in $Y$ has distance at most $2^{\delta-1}\ell$ from $Z$;
			\item for all $u,v\in Y$, if $\dist_{G[Y]}(u,v)>\delta2^\delta\ell$, then $\dist_G(u,v)>\ell$. 
		\end{itemize}
	\end{thm}
	%
	\Proof For $1\le t\le \delta$, define $f_t= 2^{\delta-t}\ell$, and $g_t = 2tf_1$. Choose a tie-breaker weighting $\lambda$ in $G$.
	Let $Z_0 = Z$. We will define a sequence of subsets $Z=Z_0\subseteq Z_1\subseteq \cdots\subseteq Z_{\delta}$ inductively, as follows. 
	For $1\le t\le \delta$, suppose that $Z_0\LL Z_{t-1}$ have been defined. 
	A {\em $t$-bite} is a path $P$ such that:
	\begin{itemize}
		\item $P$ is a $\lambda$-geodesic of length at least two;
		\item the ends $u_1,u_2$ of $P$ belong to $Z_{t-1}$;
		\item $P$ has length at most $f_t$; 
		\item if $t=1$ then no internal vertex of $P$ is in $Z$; and
		\item if $t\ge 2$ then $\dist_{Z_{t-1}}(u_1,u_2)> g_t$.
	\end{itemize}
	Let $Z_{t}$ be the union of $Z_{t-1}$ and the vertex sets of all $t$-bites. This completes the inductive definition of 
	$Z_0\LL Z_{\delta}$.
	We observe:
	\\
	\\
	(1) {\em For $1\le t\le \delta$, if $u,v\in Z_{t-1}$, and $\dist_G(u,v)\le f_t$, and either $t = 1$ or $\dist_{Z_{t-1}}(u,v)> g_t$, 
		then $\dist_{Z_t}(u,v)=\dist_G(u,v)\le f_t$.}
	\\
	\\
	To see this, we may assume that $u,v$ are distinct and nonadjacent. Let $P$ be the $\lambda$-geodesic between $u,v$.
	If $t=1$, then every vertex of $P$ not in $Z$ belongs to a $1$-bite and hence to $Z_1$, and so $V(P)\subseteq Z_t$. If $t\ge 2$,
	then $P$ is a $t$-bite, so again $V(P)\subseteq Z_t$. Hence, in either case, 
	$$\dist_{Z_t}(u,v)=|E(P)|=\dist_G(u,v)\le f_t.$$ 
	This proves (1).
	
	\bigskip
	
	For $1\le t\le \delta$ and for each $v\in Z_{t}\setminus Z_{t-1}$, $v$ belongs to the interior of a $t$-bite (perhaps of several
	$t$-bites). Choose one such $t$-bite, arbitrarily, and denote it by $S_v$. Thus, $S_v$ is defined for each $v\in Z_{\delta}\setminus Z$.
	
	For each $v\in Z$, let $L_v$ be the one-vertex subgraph with vertex $v$. 
	Inductively, for $1\le t\le \delta$, 
	if $v\in Z_{t}\setminus Z_{t-1}$, we define $L_v = S_v\cup L_{u_1}\cup L_{u_2}$, where $u_1,u_2$ are 
	the ends of $S_v$. (These graphs $L_v$ will turn out to be $Z$-leaved subdivisions of the binary tree $H_t$; so that will show that $t$ cannot be large.)
	\\
	\\
	(2) {\em For $1\le t\le \delta$, if $v\in Z_{t}$, then every vertex $w$ of $L_v$ is joined to $v$ by a path of $L_v$
		with length at most $f_t+f_{t-1}+\cdots+f_1=2f_1-f_t$.}
	\\
	\\
	We proceed by induction on $t$. Let $w\in V(L_v)$. If $w\in V(S_v)$, the statement is clear, since $S_v$ has length at most $f_t$,
	so we may assume that $w\in V(L_u)$ for some end $u$ of $P$. Since $u\in Z_{t-1}$, it follows inductively that 
	$w$ is joined to $u$ by a path of $L_u$ of length at most  $f_{t-1}+f_{t-2}+\cdots+f_1$. Since $L_u$ is a subgraph of $L_v$, the claim follows by adjoining to this path
	the subpath of $S_v$ between $v,u$. This proves (2).
	\\
	\\
	(3) {\em For $1\le t\le \delta$ and every $t$-bite $P$ with ends $u_1,u_2$, $L_{u_1}$ is vertex-disjoint from $L_{u_2}$.}
	\\
	\\
	We may assume that $t\ge 2$. Suppose that there exists $w\in V(L_{u_1}\cap L_{u_2})$. By (2), 
	$$\dist_{Z_{t-1}}(w,u_i)\le \dist_{L_{u_i}}(w,u_i)\le 2f_1$$
	for $ i = 1,2$, and so $\dist_{Z_{t-1}}(u_1,u_2)\le 4f_1$.
	Since $\dist_{Z_{t-1}}(u_1,u_2)>g_t\ge 4f_1$ (because $t\ge 2$)
	there is no such vertex $w$. This proves (3).
	\\
	\\
	(4) {\em For $2\le t\le \delta$, if $P$ is a $t$-bite with ends $u_1,u_2$, then $u_1,u_2\in Z_{t-1}\setminus Z_{t-2}$.}
	\\
	\\
	Suppose first that 
	$u_1,u_2\in Z_{t-2}$. Then 
	\begin{align*}
		\dist_G(u_1,u_2)&=|E(P)|\le f_t\le f_{t-1}, \text{ and }\\
		\dist_{Z_{t-2}}(u_1,u_2) &\ge \dist_{Z_{t-1}}(u_1,u_2)> g_t\ge g_{t-1},
	\end{align*}
	so 
	$\dist_{Z_{t-1}}(u_1,u_2)\le  f_{t-1}\le g_t$ by (1), a contradiction. 
	Thus, at least one of $u_1,u_2$, say $u_1$, belongs to $Z_{t-1}\setminus Z_{t-2}$. Suppose that $u_2\in Z_{t-2}$.
	Since $S_{u_1}$ has length at most $f_{t-1}$, it includes a subpath  
	of length at most $f_{t-1}/2$ between $u_1$ and one of the ends $w_1$ of $S_{u_1}$. It follows that $\dist_{G}(w_1, u_2)\le f_t+ f_{t-1}/2= f_{t-1}$.
	Moreover, 
	$$\dist_{Z_{t-1}}(u_1,u_2)\le \dist_{Z_{t-1}}(w_1,u_1) + \dist_{Z_{t-1}}(w_1,u_2),$$
	and so 
	$$\dist_{Z_{t-2}}(w_1,u_2)\ge \dist_{Z_{t-1}}(w_1,u_2)\ge \dist_{Z_{t-1}}(u_1,u_2) -\dist_{Z_{t-1}}(w_1,u_1)> g_t-f_{t-1}/2\ge g_{t-1}.$$
	By (1), $\dist_{Z_{t-1}}(w_1,u_2)\le f_{t-1}$. 
	It follows that $\dist_{Z_{t-1}}(u_1,u_2)\le  f_{t-1} + f_{t-1}/2\le g_t$,
	contradicting that  $\dist_{Z_{t-1}}(u_1,u_2)>g_t$.
	Thus $u_2\notin Z_{t-2}$, and so $u_2\in Z_{t-1}\setminus Z_{t-2}$. 
	This proves (4).
	\\
	\\
	(5) {\em For $2\le t\le \delta$, if $P$ is a $t$-bite with ends $u_1,u_2$, then for $i = 1,2$, every vertex of $P\cap L_{u_i}$
		belongs to $Z_{t-1}\setminus Z_{t-2}$ and hence belongs to the interior of $S_{u_i}$.}
	\\
	\\
	Suppose that $w\in V(P\cap L_{u_2})$, and $w\in Z_{t-2}$. 
	By (2),
	$$\dist_{Z_{t-1}}(u_2,w)\le 2f_1-f_{t-1},$$
	and so 
	$$\dist_{Z_{t-1}}(u_1,w)>g_t-\left(2f_1-f_{t-1}\right).$$
	As in the proof of (4), there is a subpath $P_1$ of $S_{u_1}$, between $u_1$ and some $w_1\in Z_{t-2}$, of length at most
	$f_{t-1}/2$. Hence 
	\begin{align*}
		\dist_{Z_{t-2}}(w_1,w)&\ge \dist_{Z_{t-1}}(w_1,w)> g_t-(2f_1-f_{t-1}) - f_{t-1}/2\ge g_{t-1}, \text{ and}\\
		\dist_G(w_1,w)&\le |E(P)|+|E(P_1)|\le f_t+f_{t-1}/2= f_{t-1}.
	\end{align*}
	By (1), $\dist_{Z_{t-1}}(w_1,w)\le  f_{t-1}$,
	and so 
	\begin{align*}
		\dist_{Z_{t-1}}(u_1,u_2)&\le  \dist_{Z_{t-1}}(u_1,w_1)+\dist_{Z_{t-1}}(w_1,w)+\dist_{Z_{t-1}}(w,u_2)\\
		&\le f_{t-1}/2+f_{t-1}+\left(2f_1-f_{t-1}\right)\le g_t,
	\end{align*}
	a contradiction. This proves (5).
	\\
	\\
	(6) {\em For $1\le t\le \delta$, if $v\in Z_t\setminus Z_{t-1}$, then $L_v$ (rooted at $v$) is a $(<f_1)$-subdivision of the binary tree $H_{t}$, and all its leaves are in $Z$.}
	\\
	\\
	We proceed by induction on $t$. If $t=1$ the statement is clear, so we assume that $t\ge 2$.
	Let $S_v$ have ends $u_1,u_2\in Z_{t-1}$.
	By (4),  $u_1,u_2\in Z_{t-1}\setminus Z_{t-2}$.
	From the inductive hypothesis, for $i = 1,2$, 
	$L_{u_i}$ (rooted at $u_i$) is a $(<f_{1})$-subdivision of the binary tree $H_{t-1}$, and all its leaves are in $Z$.
	By (3), $L_{u_1}, L_{u_2}$ are vertex-disjoint.
	By (5), for $i = 1,2$, every vertex of $S_v$ in $L_{u_i}$ belongs to $Z_{t-1}\setminus Z_{t-2}$ and hence belongs to the interior of $S_{u_i}$. But 
	since $S_{u_i}$ and $S_v$ are both $\lambda$-geodesics, their intersection is also a $\lambda$-geodesic $A_i$ say. Thus, for $i = 1,2$,
	$A_i$ is a 
	subpath of $S_v$ containing the end $u_i$ of $S_v$, possibly of length zero;
	and $A_1$ is vertex-disjoint from $A_2$ since $v\in V(S_v)$ and $v\notin V(S_{u_i})$ for $i = 1,2$.
	Let $A_i$ have ends $u_i,a_i$ for $i = 1,2$, and let $A_0$ be the subpath of $S_v$ between $a_1,a_2$. Thus, $S_v$ is the concatenation of the paths $A_1,A_0,A_2$. Moreover, for $i = 1,2$, every vertex of $A_i$ belongs to 
	$Z_{t-1}\setminus Z_{t-2}$ by (5); so $A_i$ contains neither end of $S_{u_i}$. 
	The interior of the path $A_0$ is disjoint from $L_{u_1}\cup L_{u_2}$ and contains $v$, and so this proves (6).
	
	\bigskip
	
	We want to show that each $L_v$ is $Z$-leaved, and to do so it only remains to prove that $L_v$ has no vertices in $Z$ except its leaves. For that we use the
	following.
	\\
	\\
	(7) {\em For $1\le t\le \delta$, if $v\in Z_{t}$ and $v'\in V(S_v)$, then $\dist_{L_v}(v',Z)\le (f_t+f_{t-1}+\cdots+f_1)/2=f_1-f_t/2$.}
	\\
	\\
	We proceed by induction on $t$. Since $S_v$ has length at most $f_t$, $v'$ is joined to an end $u$ of $S_v$ by a subpath of $S_v$
	with length at most $f_t/2$. If $t=1$ then $u_1\in Z$ and the claim holds, so we assume that $t\ge 2$. Inductively,
	$\dist_{L_{u}}(u,Z)\le (f_{t-1}+\cdots+f_1)/2$, and the claim follows. This proves (7).
	\\
	\\
	(8) {\em For $1\le t\le \delta$, if $v\in Z_t\setminus Z_{t-1}$, then the tree $L_v$ is $Z$-leaved.}
	\\
	\\
	We proceed by induction on $t$. The result is true if $t=1$ from the definition of a 1-bite, so we assume that $t\ge 2$.
	Let $u_1,u_2$ be the ends of $S_v$. 
	Inductively, the trees $L_{u_1}, L_{u_2}$ are $Z$-leaved, so we just need to show that no internal vertices of $S_v$ belong to $Z$.
	
	Suppose that $z\in Z$ is 
	an internal vertex of $S_v$. As in (6), $S_v$ is the concatenation of a subpath $A_1$ of $S_{u_1}$ between $u_1$ and some $a_1$, 
	a path $A_0$ say between $a_1,a_2$, and a subpath $A_2$ of $S_{u_2}$ between $a_2,u_2$, where:
	\begin{itemize}
		\item $a_i$ is an internal vertex of $S_{u_i}$ for $i = 1, 2$; and
		\item no internal vertex of $A_0$ belongs to $V(L_{u_1})\cup V(L_{u_2})$. 
	\end{itemize}
	It follows that $z,v$ are both internal vertices of $A_0$. Let $B_i$ be the subpath of $A_0$ between $z,a_i$ for $i = 1,2$. 
	By (7), for $i = 1, 2$ there exists $z_i\in V(L_{u_i})\cap Z$ and a path $Q_i$ of $L_{u_i}$ between $a_i,z_i$ of length at most $f_1-f_{t-1}/2$.
	Hence $\dist_G(z,z_i)\le (f_1-f_{t-1}/2)+f_t=f_1$ since $B_i$ has length at most $f_t$,
	and consequently the $\lambda$-geodesic between $z,z_1$ has length at most $f_1$. Each vertex of this $\lambda$-geodesic either belongs to $Z$ or to the interior of a $1$-bite, and therefore belongs to $Z_1$; and so 
	$\dist_{Z_1}(z_i,z)\le f_1$ for $i = 1,2$. Consequently,
	$$\dist_{Z_{t-1}}(z_1,z_2)\le \dist_{Z_1}(z_1,z_2)\le 2f_1.$$ 
	Now, for $i = 1,2$, 
	$$\dist_{Z_{t-1}}(z_i,u_i)\le |E(Q_i)|+|E(A_i)|\le (f_1-f_{t-1}/2)+f_t=f_1,$$
	since $A_i$ is a subpath of $S_v$ and so has length at most $f_t$. Hence $\dist_{Z_{t-1}}(u_1,u_2)\le 4f_1\le g_t$ (since $t\ge 2$),
	a contradiction. This proves (8).

	\bigskip
	
	Since, by hypothesis, there is no $Z$-leaved $(<f_1)$-subdivision of $H_{\delta}$ in $G$, it follows from (8) that if $Z_t\setminus Z_{t-1}\ne \emptyset$ then $t\le \delta-1$.
	In particular, $Z_{\delta} = Z_{\delta-1}$, and it follows that every $\delta$-bite is contained in $G[Z_{\delta-1}]$. We deduce that 
	if $u,v\in Z_{\delta-1}$ are distinct, and $\dist_G(u,v)\le f_\delta$, then $\dist_{Z_{\delta-1}}(u,v)\le g_{\delta}$. Moreover, by (7), every vertex 
	$v\in Z_{\delta-1}$ has distance at most $f_1$ from $Z$, since $L_v$ has a vertex in $Z$
	by (6). Thus, setting $Y=Z_{\delta-1}$
	satisfies the theorem. This proves \ref{newkey}.~\bbox
	
	In \ref{newkey} we assumed that no subdivision of $H_\delta$ has all its leaves in $Z$ and no other vertices in $Z$, but now we need to 
	refine that. 
	
	\begin{thm}\label{settle}
		Let $d\ge 1$ be an integer, let $G$ be a graph, let $\Gamma$ be a subgraph of $G$, such that no subgraph is a $V(\Gamma)$-leaved subdivision of 
		$H_{d}$ whose leaves all belong to different components of $\Gamma$.
		Suppose that 
		$M_1,M_2\LL M_t$
		are paths of $G$, each of length at most $\ell$. Let  $\Gamma_i=\Gamma\cup (M_1\cupcup M_i)$ for $0\le i\le t$; and suppose in addition that
		for each $i\ge 1$, the ends of $M_i$ lie in different components  of
		$\Gamma_{i-1}$, and none of its internal vertices lie in $V(\Gamma_{i-1})$. 
		Then for each $v\in V(\Gamma_t)\setminus V(\Gamma)$, either $v$ lies in the interior 
		of some $M_i$ with both ends in $V(\Gamma)$, or there are at least
		three components $C$ of $\Gamma$ such that $v$ is joined to $C$ by a path in $\Gamma_{t}$ of length at most $(d+1)(\ell-1)$.
	\end{thm}
	\Proof For $h\ge1$, we say that a $V(\Gamma)$-leaved subdivision of $H_h$ in $G$ is {\em $\Gamma$-pruned} if its leaves are all in distinct components of $\Gamma$. We say the {\em height} of each vertex in $V(\Gamma)$ is zero; and inductively, 
	for $1\le i\le t$, let us say that for each vertex in the interior of $M_i$, its {\em height} is one more than the minimum of
	the heights of $u_1,u_2$, where $u_1,u_2$ are the ends of $M_i$. Then:
	\\
	\\
	(1) {\em For each $i\ge 0$ and each $v\in V(\Gamma_i)$ with height at least $h\ge 1$, there is
		a subgraph of $\Gamma_i$ that is a $\Gamma$-pruned subdivision 
		of $H_{h}$ rooted at $v$.
	}
	\\
	\\
	We use induction on $h$. The statement is clear if $h=1$, so we assume $h\ge 2$. We may assume that $i$ is minimum such that 
	$v\in V(\Gamma_i)$,
	and consequently $v$ belongs to the interior of $M_i$.  Let $u_1,u_2$ be the ends of $M_i$, joining components $C_1,C_2$ of
	$\Gamma_{i-1}$. Thus, $u_1,u_2$ have height at least $h-1\ge 1$.
	From the inductive hypothesis there is a subgraph $L_j$ of $C_j$ rooted at $u_j$
	that is a $\Gamma$-pruned subdivision of $H_{h-1}$.
	But $L_1,L_2$ are disjoint, and every component of $\Gamma$ containing a leaf of $L_1$ is different from every component of $\Gamma$ containing a leaf of $L_2$, since they belong to different components 
	of $\Gamma_{i-1}$. Moreover, $L_1,L_2$ are both disjoint from the interior of $M_i$, since the latter is disjoint from $V(\Gamma_{i-1})$. 
	Consequently, $L_1\cup L_2\cup M_i$ (rooted at $v$) is the desired 
	subdivision    
	of $H_{h}$. This proves (1). 
	
	\bigskip
	From (1) and the hypothesis, it follows that every vertex has height at most $d-1$.
	\\
	\\
	(2) {\em For each $i\ge 0$ and each $v\in V(\Gamma_i)$ with height $h\ge 0$, 
		$v$ is joined to $V(\Gamma)$ by a path in $\Gamma_i$ of length at most $h(\ell-1)$.}
	\\
	\\
	We prove this by induction on $h\ge 0$. If $h=0$, the statement is clear, so we assume that $h\ge 1$.
	Choose $i$ minimum with $v\in V(\Gamma_i)$. 
	Then  $v$ is 
	joined to a vertex $u$ of height $h-1$ by a path of $\Gamma_i$ of length 
	at most $\ell-1$ (a subpath of $M_i$); and from the inductive hypothesis, $u$ is joined to $V(\Gamma)$ by a path in $\Gamma_{i-1}$
	(and hence of $\Gamma_i$)
	of length at most $(h-1)(\ell-1)$. Consequently $v$ is joined to $V(\Gamma)$ by a path in $\Gamma_i$
	of length at most $h(\ell-1)$. This proves (2). 
	\\
	\\
	(3) {\em  For each $i\ge 1$ and each $v\in V(\Gamma_i)\setminus V(\Gamma)$, there are at least two components $C$ of $\Gamma$
		such that $v$ 
		is joined to $C$ by a path in $\Gamma_i$ of length at most $d(\ell-1)$.}
	\\
	\\
	Choose $i$ minimum with $v\in V(\Gamma_i)$. Thus, $v$ belongs to the interior of $M_i$;
	let $M_i$ have ends $u_1,u_2$. Each of $u_1,u_2$ belongs to $V(\Gamma_{i-1})$, and since $u_1,u_2$ have height at most $d-1$ by (1), it 
	follows from (2) that each of $u_1,u_2$
	is joined to a component of $\Gamma$ by a path in $\Gamma_{i-1}$ of length at most  $(d-1)(\ell-1)$. These components of $\Gamma$
	are different since $u_1,u_2$ are in different components of $\Gamma_{i-1}$, and $v$ is joined to each of them by 
	a path in $\Gamma_i$ of length at most $d(\ell-1)$. This proves (3).

	\bigskip
	
	In particular, for each $v\in V(M_1\cupcup M_t)\setminus V(\Gamma)$, $v$ has height at least one; choose $i$ minimum with $v\in V(\Gamma_i)$. 
	Thus, $v$ belongs to the interior of $M_i$;
	let $M_i$ have ends $u_1,u_2$. If $u_1,u_2$ both have height zero then $M_i$ has both 
	ends in $V(\Gamma)$ and the theorem holds; so we assume that $u_1\notin V(\Gamma)$. By (3), 
	there are at least two components $C_1,C_2$ of $\Gamma$ that are joined to 
	$u_1$
	by a path in $\Gamma_{i-1}$ of length at most $d(\ell-1)$; and by (2), there is a component $C_3$ of $\Gamma$ such that 
	$u_2$
	is joined to $C_3$ by a path in $\Gamma_{i-1}$ of length at most $(d-1)(\ell-1)$. Moreover, $C_3\ne C_1,C_2$ since $u_1,u_2$
	belong to different components of $\Gamma_{i-1}$. Consequently there are at least 
	three components $C$ of $\Gamma$ such that $v$ is joined to $C$ by a path in $\Gamma_t$ of length at most $(d+1)(\ell-1)$.
	This proves \ref{settle}.~\bbox

	
	\section{Augmenting paths}
	
	Let us extend the definition
	of $\dist_G(u,v)$  a little, to accommodate vertices $u,v\notin V(G)$: if one of $u,v\notin V(G)$ then $\dist_G(u,v)=\infty$.
	
	Some more notation: if $P$ is a path and $u,v\in V(P)$, $P[u,v]$ denotes the subpath between $u,v$. If $\mac P$ is a set of vertex-disjoint paths of a graph $G$, we denote $P_1\cupcup  P_k$ by $U\mac P$, and its vertex set by $V\mac P$.
	Let $G$ be a graph, let $S,T\subseteq V(G)$ be disjoint, and let $\mac P=\{P_1\LL P_k\}$ be a set of $k$ vertex-disjoint $S$-$T$ paths, with 
	$$V\mac P\cup S\cup T=V(G),$$ 
	such that for $1\le h\le k$, no proper subpath of $P_h$ is an $S$-$T$ path.  Let $P_h$ have ends $s_h\in S$ and $t_h\in T$.
	If $u,v\in V(P_h)$ are distinct, and  $v$ belongs to $P_h[u,t_h]$, we say that $v$ is {\em later than $u$ in $P_h$}, and 
	{\em $u$ is earlier than $v$ in $P_h$}.
	
	It is an elementary theorem (a special case of the theory of augmenting paths)
	that: 
	\begin{thm}\label{disjtpaths}
		Given $G,S,T$ and $\mac P = \{P_1\LL P_k\}$ as above,
		the following are equivalent:
		\begin{enumerate}
			\item For every choice of $v_i\in V(P_i)$ for $1\le i\le k$, there is an edge $ab$ of $G$ with $a,b\notin 
			\{v_1\LL v_k\}$,  such that
			\begin{itemize}
				\item  either $a\in S\setminus V\mac P$ or for some $h\in \{1\LL k\}$, $a\in V(P_h)$,
				and $a$ is earlier than $v_h$ in $P_h$, and 
				\item either
				$b\in T\setminus V\mac P$ or for some $h\in \{1\LL k\}$, $b\in V(P_h)$ and $b$ is later than $v_h$ in $P_h$.
			\end{itemize}
			
			\item There is a sequence $a_1b_1,a_2b_2\LL a_nb_n$ of oriented edges of $G$, not in $E(P_1\cupcup P_k)$,  such that 
			\begin{itemize}
				\item $a_1\in S\setminus V\mac P$, and $b_n\in T\setminus V\mac P$;
				\item for $1\le i<n$, $b_i, a_{i+1}$ belong to the same path $P_h$ say (where $1\le h\le k$), and 
				$a_{i+1}$ is earlier than $b_i$ in $P_h$.
			\end{itemize}
			
			\item There is a sequence $a_1b_1,a_2b_2\LL a_nb_n$ as above, satisfying in addition that for $1\le h\le k$, and $1\le i<j\le n$,
			if $u\in \{a_i,b_i\}\cap V(P_h)$ and $v\in  \{a_j,b_j\}\cap V(P_h)$, then either 
			\begin{itemize}
				\item $u$ is earlier than $v$ in $P_h$, or
				\item $b_i=u=v=a_j$; or 
				\item $b_i=u$ and $a_j=v$ and $j=i+1$. 
			\end{itemize}
			\item There are $k+1$ vertex-disjoint $S$-$T$ paths in $G$.
		\end{enumerate}
	\end{thm}
	We do not actually need this theorem, and we mention it just for comparison with the more complicated results that we will need.
	
	Let $S,T$ be sets, and let $\mac P = \{P_1\LL P_k\}$ be a set of $k$ vertex-disjoint $S$-$T$ paths, each a minimal $S$-$T$ path.
	(Possibly some $P_i$ has length zero, and then its vertex is in $S\cap T$.) We call $(S,T,\mac P)$ a {\em setting}. 
	Let $F_0$ be the set of all ordered pairs of vertices $ab$ with   
	$a,b\in V\mac P\cup S\cup T$. 
	
	Let us fix some setting $(S,T,\mac P)$ where $\mac P = \{P_1\LL P_k\}$.
	Let $c\ge 0$ be an integer. A {\em $c$-barrier} (in the setting) is a $k$-tuple $Q_1\LL Q_k$, where $Q_h$ is a subpath of $P_h$ of length at most $c$. 
	We say $ab\in F_0$ {\em jumps} a $c$-barrier $Q_1\LL Q_k$ (in the setting) if $a,b\notin V(Q_1\cupcup Q_k)$, and 
	\begin{itemize}
		\item  either $a\in S\setminus V\mac P$ or for some $h\in \{1\LL k\}$, $a\in V(P_h)$, and $a$ is earlier than each vertex of $Q_h$ in $P_h$; and
		\item either
		$b\in T\setminus V\mac P$ or for some $h\in \{1\LL k\}$, $b\in V(P_h)$ and $b$ is later than each vertex of $Q_h$ in $P_h$.
	\end{itemize}
	
	Let us say a set $F\subseteq F_0$ is {\em $c$-jumping} (in the setting $(S,T,\mac P)$) if for every $c$-barrier,
	some member of $F$ jumps the $c$-barrier.

	A {\em partial $c$-augmenting sequence to $b_n$} is a sequence $a_1b_1,a_2b_2\LL a_nb_n$ of elements of $F_0$,
	such that
	\begin{itemize}
		\item $a_1\in S\setminus V\mac P$;
		\item for $1\le i<t$, there exists $h\in \{1\LL k\}$ with $b_i, a_{i+1}\in V(P_h)$, and
		$a_{i+1}$ is earlier than $b_i$ in $P_h$, and $P_h[a_{i+1},b_i]$ has length at least $c+1$.
	\end{itemize}
	If in addition $b_n\in T\setminus V\mac P$, we call such a sequence a {\em $c$-augmenting sequence}.
	Thus, if $a_i\in S\setminus V\mac P$ and $b_i\in T\setminus V\mac P$ then $i=1=n$. 
	For $F\subseteq F_0$, the sequence is {\em in $F$} if $a_ib_i\in F$ for $1\le i\le n$.
	
	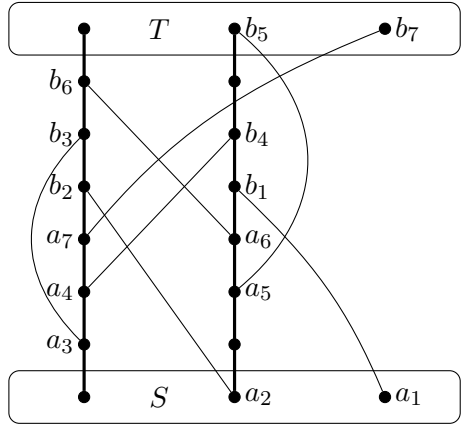
\begin{figure}[H]
		\centering
		
		\begin{tikzpicture}[yscale=.7,auto=left]
			
			\tikzstyle{every node}=[inner sep=1.5pt, fill=black,circle,draw]
			
			\draw[rounded corners] (-1,-.5) rectangle (5,.5);
			\draw[rounded corners] (-1,6.5) rectangle (5,7.5);
			\node (u0) at (0,0) {};
			\node (v0) at (2,0) {};
			\node (u1) at (0,1) {};
			\node (v1) at (2,1) {};
			\node (u2) at (0,2) {};
			\node (v2) at (2,2) {};
			\node (u3) at (0,3) {};
			\node (v3) at (2,3) {};
			\node (u4) at (0,4) {};
			\node (v4) at (2,4) {};
			\node (u5) at (0,5) {};
			\node (v5) at (2,5) {};
			\node (u6) at (0,6) {};
			\node (v6) at (2,6) {};
			\node (u7) at (0,7) {};
			\node (v7) at (2,7) {};

			\node (a1) at (4,0) {};
			\node (b7) at (4,7) {};
			
			\draw[very thick] (u0)--(u7);
			\draw[very thick] (v0)--(v7);
			
			\foreach \from/\to in {v0/u4,u2/v5,v3/u6}
			\draw [-] (\from) -- (\to);
			\draw[-] (u1) to [bend left = 35] (u5);
			\draw[-] (v2) to [bend right = 40] (v7);
			
			\draw[-] (u3) to [bend left = 15] (b7);
			\draw[-] (a1) to [bend right = 10] (v4);

			\tikzstyle{every node}=[]
			\draw[right] (a1) node []           {$a_1$};
			\draw[right] (v4) node []           {$b_1$};
			\draw[right] (v0) node []           {$a_2$};
			\draw[left] (u4) node []           {$b_2$};
			\draw[left] (u1) node []           {$a_3$};
			\draw[left] (u5) node []           {$b_3$};
			\draw[left] (u2) node []           {$a_4$};
			\draw[right] (v5) node []           {$b_4$};
			\draw[right] (v2) node []           {$a_5$};
			\draw[right] (v7) node []           {$b_5$};
			\draw[right] (v3) node []           {$a_6$};
			\draw[left] (u6) node []           {$b_6$};
			\draw[left] (u3) node []           {$a_7$};
			\draw[right] (b7) node []           {$b_7$};
			
			\draw[] node at (1,0)            {$S$};
			\draw[] node at (1,7) []           {$T$};

		\end{tikzpicture}
		
		\caption{$P_1,P_2$ are the two paths of thick edges. With $k=2$, the sequence $a_1b_1\LL a_7b_7$ is minimal 2-augmenting, but it is not 1-separated. For every choice of three 
			vertex-disjoint  $S$-$T$ paths, some edge of $P_1\cup P_2$ joins two of them.} \label{fig:aug}
	\end{figure}
	
	We begin with:
	\begin{thm}\label{organize}
		Let $(S,T,\mac P)$ be a setting, with $\mac P = \{P_1\LL P_k\}$, and let $c\ge 0$ be an integer. With $F_0$ as before, let $F\subseteq F_0$. Then the following are equivalent:
		\begin{itemize}
			\item $F$ is $c$-jumping;
			\item there is a $c$-augmenting sequence of elements of $F$.
		\end{itemize}
	\end{thm}
	\Proof
	We show first that the second statement implies the first. To see this, assume that
	the sequence $a_1b_1,a_2b_2\LL a_nb_n$ of pairs in $F$ is $c$-augmenting, 
	and let $Q_1\LL Q_k$ be a $c$-barrier. Choose $i$ maximum such that either $a_i\in S\setminus V\mac P$, or for some 
	$h\in \{1\LL k\}$, $a_i\in V(P_h)\setminus V(Q_h)$ and $a_i$ is earlier in $P_h$ than each vertex of $Q_h$. 
	If $b_i\in T\setminus V\mac P$
	then $a_ib_i$ jumps the $c$-barrier, so we assume that $b_i\in V(P_j)$ for some $j\in \{1\LL k\}$. Consequently $i<n$, and
	$a_{i+1}\in V(P_j)$, earlier than $b_i$ in $P_j$. From the maximality of $i$, there exists
	$q\in V(Q_j)$ such that  $a_{i+1}$ is not earlier than $q$ in $P_j$.
	Since $P_j[a_{i+1}, b_i]$ has length at least $c+1$, it follows that $b_i$ is later than $q$ in $P_j$,
	and $P_j[q, b_i]$ has length at least $c+1$. Since $Q_j$ has length at most $c$, it follows that $b_i$
	is later in $P_j$ than every vertex of $Q_j$; and so $a_ib_i$ jumps the $c$-barrier. This proves the first statement.
	
	To show the converse, suppose that $F$ is $c$-jumping, and for $1\le h\le k$, choose $v_h\in V(P_h)$ with $P_h[s_h,v_h]$ maximal such that either $v_h=s_h$ or there is a partial $c$-augmenting sequence to $v_h$ in $F$.
	For $1\le h\le k$, let $Q_h$ be the maximal subpath of $P_h[s_h,v_h]$ with
	length at most $c$, such that one of its ends is $v_h$. Thus, $Q_1\LL Q_k$ is a $c$-barrier, and so, since $F$ is $c$-jumping, 
	some $ab\in F$ jumps
	this $c$-barrier.  Suppose first that $a\in S\setminus V\mac P$. If $b\in T\setminus  V\mac P$, then $ab$ is
	a $c$-augmenting sequence, so we assume that $b\in V(P_h)$ for some $h\in \{1\LL k\}$. Since $ab$ jumps the $c$-barrier, it follows that
	$b$ is later than $v_h$ in $P_h$, contradicting the choice of $v_h$, since $ab$ is a partial $c$-augmenting sequence to $b$.
	Thus, we may assume that  for some $h\in \{1\LL k\}$, $a\in V(P_h)$, and $a$ is earlier than each vertex of $Q_h$ in $P_h$.
	Since $a\notin V(Q_h)$,
	it follows from the maximality of $Q_h$ that $Q_h$ has length exactly $c$, and therefore $P_h[a,v_h]$ has length at least $c+1$. 
	Let
	$a_1b_1\LL a_sb_s$ be a partial $c$-augmenting sequence to $v_h$ in $F$.
	Consequently $a_1b_1\LL a_sb_s,ab$
	is a partial $c$-augmenting sequence to $b$ in $F$.
	If $b\notin T\setminus V\mac P$, then, since $ab$ jumps the $c$-barrier, there exists $h'\in \{1\LL k\}$ such that
	$b\in P_{h'}[v_{h'},t_{h'}]$ and $b\ne v_{h'}$; but this contradicts the definition of $v_{h'}$. Thus,
	$b\in T\setminus V\mac P$, and so $a_1b_1\LL a_sb_s,ab$
	is a $c$-augmenting sequence in $F$. This proves \ref{organize}.~\bbox
	
	This provides an analogue of the first two bullets of \ref{disjtpaths}, and the next result gives an analogue of the third bullet.
	\begin{thm}\label{minimal}
		Let $(S,T,\mac P)$ be a setting, with $\mac P = \{P_1\LL P_k\}$, and let $c\ge 0$ be an integer. Let $F\subseteq F_0$ be  
		$c$-jumping, and choose a $c$-augmenting sequence $a_1b_1\LL a_nb_n$ of elements of $F$, with $n$ minimum. 
		For $1\le h\le k$, and $1\le i<j\le n$,
		if $u\in \{a_i,b_i\}\cap V(P_h)$ and $v\in  \{a_j,b_j\}\cap V(P_h)$, then either 
		\begin{itemize}
			\item $u$ is earlier than $v$ in $P_h$; or
			\item $b_i=u$ and $v=a_j$ and $P_h[u,v]$ has length at most $c$; or
			\item $b_i=u$ and $v=a_j$ and $j=i+1$.
		\end{itemize}
	\end{thm}
	\Proof
	Suppose that
	$1\le h\le k$, and $1\le i<j\le n$, and $u\in \{a_i,b_i\}\cap V(P_h)$ and $v\in  \{a_j,b_j\}\cap V(P_h)$, and $u$ is not earlier than $v$ in $P_h$.
	If $u=a_i$ and $v=a_j$, then $i\ge 2$ and
	$$a_1b_1\LL a_{i-1}b_{i-1}, a_jb_j\LL a_nb_n$$
	is a $c$-augmenting sequence in $F$, contrary to the minimality of $n$. Similarly, if $u=b_i$ and $v=b_j$, then
	$$a_1b_1\LL a_{i}b_{i}, a_{j+1}b_{j+1}\LL a_nb_n$$
	is a $c$-augmenting sequence, a contradiction; and if $u=a_i$ and $v=b_j$, then $i\ge 2$ and $j\le n-1$ and 
	$$a_1b_1\LL a_{i-1}b_{i-1}, a_{j+1}b_{j+1}\LL a_nb_n$$
	is a $c$-augmenting sequence, a contradiction. Thus, we assume that $u=b_i$ and $v=a_j$. If $P_h[u,v]$
	has length at least $c+1$, then 
	$$a_1b_1\LL a_{i}b_{i}, a_jb_j\LL a_nb_n$$
	is a $c$-augmenting sequence, and so $j=i+1$; and otherwise $P_h[u,v]$ has length at most $c$. In either case the result holds. This proves \ref{minimal}.~\bbox
	
	The results \ref{organize} and \ref{minimal} do not quite provide an analogue of \ref{disjtpaths}, because we have no 
	counterpart to the fourth statement of  \ref{disjtpaths}, the existence of $k+1$ vertex-disjoint $S$-$T$ paths. One might hope that 
	\begin{itemize}
		\item {\em In the graph obtained from $U\mac P$ by adding the remainder of $S\cup T$ as extra vertices and the pairs in $F$ as edges,
			there exist $k+1$ $S$-$T$ paths, such that no two of them are joined by a path of $U\mac P$ of length at most $c$.}
	\end{itemize}
	could be added to the list of equivalent statements given by \ref{organize} and \ref{minimal} to give an analogue of the fourth
	statement of \ref{disjtpaths}, but that is wrong. 
	This statement does imply the statements of \ref{organize}, but the 
	reverse implication does not hold. For instance, the graph of Figure \ref{fig:aug} with $k=2$ gives a 2-augmenting sequence,
	and yet for every three vertex-disjoint $S$-$T$ paths, some two of them are joined by one of the edges of $P_1\cup P_2$, which is more than we needed for a counterexample.

	The property given by \ref{minimal} implies that for each $h\in \{1\LL k\}$, the vertices $a_i$ that lie in $P_h$ 
	are all distinct and in order in $P_h$,
	but it does not imply that they are far apart in $P_h$. For instance, if $k=1$ and $P_1$ has vertices
	$s_1=v_1\CC v_n=t_1$, and $S=\{s_1,s_2\}$
	and $T=\{t_1,t_2\}$, and $F$ is the union of $\{s_2v_{c+2}, v_{n-c-1}t_2\}$ and
	the pairs  $v_iv_{i+c+2}$ for $1\le i\le n-c-2$, then the only $c$-augmenting sequence in $F$ uses all of $F$. Nevertheless, we can 
	arrange that the $a_i$'s are far apart, and the $b_j$'s are far apart, by sacrificing some of the jumping power. We show this in 
	two steps: first we arrange that the $b_j$'s are far apart, in the following.

	We recall that $\dist_{U\mac P}(b,b')=\infty$ unless $b,b'\in V\mac P$ and $b,b'$ belong to the same component of $U\mac P$.
	\begin{thm}\label{separatetops}
		Let $p,q\ge 0$ be integers, and let $F\subseteq F_0$ be $(p+q)$-jumping.
		Then there exists $D\subseteq F$ that is $p$-jumping, such that 
		if $ab,a'b'\in D$ are distinct then $\dist_{U\mac P}(b,b')>q$. 
	\end{thm}
	\Proof We will use a modified version of the second half of the proof of \ref{organize}.
	We say a partial 
	$p$-augmenting sequence $a_1b_1\LL a_sb_s$ is {\em end-separated} if 
	$\dist_{U\mac P}(b_i,b_j)>q$ for all distinct $i,j\in \{1\LL s\}$.
	By \ref{organize} it suffices to show that there is an end-separated $p$-augmenting sequence in $F$.
	
	For $1\le h\le k$, choose $v_h\in V(P_h)$ with $P_h[s_h,v_h]$ maximal such that either $v_h=s_h$ or there is an end-separated partial 
	$p$-augmenting sequence to $v_h$ in $F$.
	For $1\le h\le k$, let $Q_h$ be the maximal subpath of $P_h$ containing $v_h$, such that $Q_h\cap P_h[s_h, v_h]$ has 
	length at most $p$, and $Q_h\cap P_h[v_h,t_h]$ has length at most $q$.
	Thus, $Q_1\LL Q_k$ is a $(p+q)$-barrier, and so, since $F$ is $(p+q)$-jumping, some $ab\in F$ jumps this $(p+q)$-barrier.  
	Suppose first that $a\in S\setminus V\mac P$. If $b\in T\setminus  V\mac P$, then $ab$ is
	an end-separated $p$-augmenting sequence, so we assume that $b\in V(P_h)$ for some $h\in \{1\LL k\}$. Since $ab$ jumps the $(p+q)$-barrier, it follows that
	$b$ is later than $v_h$ in $P_h$, contradicting the choice of $v_h$, since $ab$ is an end-separated partial $p$-augmenting sequence to $b$ in $F$.
	
	Thus, we may assume that  for some $h\in \{1\LL k\}$, $a\in V(P_h)$, and $a$ is earlier than each vertex of $Q_h$ in $P_h$.
	Let
	$a_1b_1\LL a_sb_s$ be an end-separated  partial $p$-augmenting sequence to $v_h$ in $F$. Since $a\notin V(Q_h)$,
	it follows that $Q_h\cap P_h[s_h, v_h]$ has length exactly $p$, and $P_h[a,v_h]$ has length at least $p+1$. Consequently 
	$a_1b_1\LL a_sb_s,ab$
	is a partial $p$-augmenting sequence to $b$ in $F$. 
	If $b\notin T\setminus V\mac P$, then, since $ab$ jumps the $(p+q)$-barrier, there exists $h'\in \{1\LL k\}$ such that
	$b\in P_{h'}[v_{h'},t_{h'}]$ and $b\ne  V(Q_{h'})$; but then $Q_{h'}\cap P_{h'}[v_{h'},t_{h'}]$ has length exactly $q$, and so 
	$P_{h'}[v_{h'},b]$ has length $>q$. Since each $b_i$ in $V(P_{h'})$ belongs to $P_{h'}[s_{h'}, v_{h'}]$ from the definition of 
	$v_{h'}$, it follows that $a_1b_1\LL a_sb_s,ab$ is an end-separated $p$-augmenting sequence to $b$ in $F$,
	contrary to the definition of $v_{h'}$. Thus,
	$b\in T\setminus V\mac P$, and so $a_1b_1\LL a_sb_s,ab$
	is an end-separated $p$-augmenting sequence in $F$.  This proves \ref{separatetops}.~\bbox
	
	Let us say a subset $D\subseteq F_0$ is {\em $\ell$-separated} if 
	$\dist_{U\mac P}(a,a')>\ell$ and $\dist_{U\mac P}(b,b')>\ell$ for all
	distinct $ab,a'b'\in D$. 
	We deduce:
	\begin{thm}\label{separateall}
		In the same notation, let $c\ge 0$ be an integer, and let $F\subseteq F_0$ be $5c$-jumping.
		Then there exists $D\subseteq F$ that is $c$-jumping and $2c$-separated.
	\end{thm}
	\Proof
	This follows from two applications of \ref{separatetops}: first, to $F$ with $(p,q) = (3c,2c)$, giving some $3c$-jumping set $F'$;
	and then to $F'$ with $S,T$ exchanged and $(p,q) = (c,2c)$. This proves \ref{separateall}.~\bbox
	
	We remark that, even if $F$ is $c$-jumping and $2c$-separated, and if $a_1b_1\LL a_nb_n$ is a minimal $c$-augmenting sequence
	of members of $F$, it is possible that $a_i=b_j$ for some pairs $i,j$.
	Even so, now we can obtain something like an analogue of the fourth statement of \ref{disjtpaths}:
	\begin{thm}\label{getpaths}
		In the same notation, let $c\ge 0$ be an integer, and let $F\subseteq F_0$ be $c$-jumping and $2c$-separated.
		Let $H$ be obtained from $U\mac P$ by adding the remainder of $S\cup T$ as vertices, and the pairs in $F$ as edges.
		Then there exist $k+1$ vertex-disjoint $S$-$T$ paths in $H$,
		such that no two of them are joined by a path of $U\mac P$ of length at most $c$.
	\end{thm}
	\Proof
	By \ref{minimal}, there is a $c$-augmenting sequence $a_1b_1\LL a_nb_n$ in $F$
	such that:
	\\
	\\
	(1) {\em 
		For $1\le h\le k$, and $1\le i<j\le n$,
		if $u\in \{a_i,b_i\}\cap V(P_h)$ and $v\in  \{a_j,b_j\}\cap V(P_h)$, then either
		\begin{itemize}
			\item $u$ is earlier than $v$ in $P_h$; or
			\item $b_i=u$ and $v=a_j$ and $P_h[a_j, b_i]$ has length at most $c$; or
			\item $b_i=u$ and $v=a_j$ and $j=i+1$.
		\end{itemize}
	}
	\noindent
	We deduce:
	\\
	\\
	(2) {\em Let $1\le h\le k$, and $1\le i\le n$ with $b_i\in V(P_h)$ (and hence $a_{i+1}\in V(P_h)$); then for $1\le j\le n$, if 
		$a_j$ belongs to $P_h[a_{i+1},b_i]$ then either
		$j=i+1$, or
		$j>i+1$ and $P_h[a_j,b_{i+1}]$ has length at most  $c$. Consequently there is at most one value of $j\ne i+1$ with 
		$a_j\in V(P_h[a_{i+1}, b_i])$, and any such $j$ satisfies $j\ge i+2$.
		Similarly there is at most one value of $j\ne i$ with $b_j\in V(P_h[a_{i+1}, b_i])$, and any such $j$ satisfies $j\le i-1$. 
	}
	\\
	\\
	By (1), $a_1\LL a_n$ are all distinct, and $b_1\LL b_n$ are all distinct.
	Suppose that $a_j$ belongs to  $P_h[a_{i+1},b_i]$, and $j\ne i+1$.  Thus, $a_{i+1}$ is earlier than $a_j$ in $P_h$. If $j\le i$
	then setting $u=a_j$ and $v=a_{i+1}$ in (1) yields a contradiction; so $i<j$, and hence $j\ge i+2$. 
	By (1) with $u=b_{i}$ and $v=a_j$, it follows that $P_h[a_j, b_{i}]$ has length at most $c$. Consequently, if $j'\ne j$
	also satisfies that $a_{j'}$ belongs to  $P_h[a_{i+1},b_i]$, and $j'\ne i+1$, then $P_h[a_j,a_{j'}]$
	has length at most $c$, contradicting that $F$ is $2c$-separated. This proves the first assertion of (2), and the second follows from
	the symmetry. This proves (2).
	
	\bigskip
	
	For $1\le i<n$, $b_i$ and $a_{i+1}$ both belong to the same member of $\mac P$, say $P_h$; let $R_i=P_h[a_{i+1},b_i]$. 
	\\
	\\
	(3) {\em Every vertex in $V\mac P$ belongs to at most two of $R_1\LL R_{n-1}$.}
	\\
	\\
	Suppose that some vertex $w$ of $P_h$ belongs to $R_{i}, R_{i'}, R_{i''}$, where $i<i'<i''$. Thus, 
	$$a_{i+1}, a_{i'+1},a_{i''+1}, b_i, b_{i'}, b_{i''}$$ 
	are in order in $P_h$ (and are all distinct except possibly $a_{i''+1}= b_i$), and $w\in P_h[a_{i''+1}, b_i]$. By (1) with 
	$u=b_i, v=a_{i'+1}$, $P_h[a_{i'+1}, b_i]$ has length at most $c$. But it includes $P_h[a_{i'+1},a_{i''+1}]$
	as a subpath, and this has length at least $2c+1$ since $F$ is $2c$-separated, a  contradiction. This proves (3).
	
	\bigskip
	For $i = 0,1,2$, let $X_i$ be the set of edges of $U\mac P$ that belong to exactly $i$ of $R_1\LL R_{n-1}$.
	Let $H'$ be the digraph obtained as follows:
	\begin{itemize}
		\item Direct the edges of $P_h$ from $s_h$ to $t_h$ for $1\le h\le k$, and let $H''$ be obtained from their union by adding the (directed) edges $a_ib_i$ for $1\le i\le n$.
		(Thus the only vertices of $S$ in V(H'') are $s_1\LL s_k$ and $a_1$.)  
		\item Reverse the direction of all edges in $X_2$.
		\item Delete all edges in $X_1$.
	\end{itemize}
	We claim:
	\\
	\\
	(4) {\em Every vertex of $H'$ either has outdegree one and indegree one, or has oudegree zero and indegree zero, except for $a_1, s_1\LL s_k$, which have outdegree one and indegree zero, and $t_1\LL t_k, b_n$, which have indegree one and outdegree zero.}
	\\
	\\
	Let $v\in V(H')$. The claim is true for $v$ if $v\in \{a_1,b_n\}$, so we may assume that $v\in V\mac P$.
	Let $v\in V(P_h)$ where $1\le h\le k$. 
	Suppose that $v=s_h$, and let $e$ be the edge of $P_h$ incident with $v$. 
	It follows that $v\ne b_1\LL b_n$, and either $v\ne a_1\LL a_n$ (and then $e\in X_0\subseteq E(H')$) 
	or $v=a_i$ for some $i<n$ (and then $e\in E(R_{i-1})$, and so $e\in X_1$ and $e\notin E(H')$).
	In either case $v$ has outdegree one and indegree zero in $H'$. Thus we may assume that $v\ne s_h$ and similarly
	$v\ne t_h$, and so $v$ is an internal vertex of $P_h$. Let $e_1,e_2$ be the two edges of $P_h$ 
	incident with $v$, where $e_1$ is in the subpath between $s_h$ and $v$.

	If $v\ne \{a_2\LL a_n, b_1\LL b_{n-1}\}$ then $v$ is an internal vertex of the (at most two) paths $R_i$ that contain $v$, and so
	satisfies the claim. Thus from the symmetry, we may assume that $v\in \{a_2\LL a_n\}$; let $v=a_i$ say.
	Thus  $i\ge 2$ and $v$ is an end of $R_{i-1}$, and $e_2\in E(R_{i-1})$.
	
	Suppose next that 
	$v\ne b_1\LL b_{n-1}$.
	If $e_2\in X_1$ then $e_2\notin E(H')$, and $e_1\in X_0\subseteq E(H')$, and $v$ therefore has indegree and outdegree one in $H'$;
	so we may assume that $e_2\in X_2$. Let $e_2\in E(R_j)$ say, where $1\le j\le n-1$ and $j\ne i$. Since $a_j,b_j\ne v$, it follows that
	$e_1\in E(R_j)$, and so $e_1\in X_1$ by (3) applied to $v$; but then $e_2\in E(H'), e_1\notin E(H')$, and the claim is true for $v$ 
	(because $v$ is the head of the directed edge $e_2$, since $e_2\in X_2$). 
	
	So we may assume that $v\in \{b_1\LL b_{n-1}\}$; let $v=b_j$. Hence $i\ne j$, and $v$ is an end of both $R_{i-1}, R_j$, and both $e_1,e_2\in X_1$ by (3). Hence $e_1,e_2\notin E(H')$, and again the claim holds for $v$. This proves (4). 
	
	\bigskip
	
	Let $J$ be the undirected graph underlying $H'$. From (4), each component of $J$ is either an $S$-$T$ path or a cycle or a vertex of degree zero; 
	and $a_1, s_1\LL s_k$ all belong to different 
	components. Since $a_1, s_1\LL s_k$ all have outdegree one and indegree zero in $H'$, and vice versa for $b_1, t_1\LL t_k$, it follows that 
	there are $k+1$ vertex-disjoint $S$-$T$ paths $P_1'\LL P'_{k+1}$ in $J$, each a component of $J$. It remains to show that no two of these paths 
	are joined by a path  
	of $U\mac P$ with length at most $c$. Suppose that $Q$ is such a path; and we can assume that no internal vertex
	of $Q$ belongs to any of $P_1'\LL P'_{k+1}$. Consequently the first and last edges of $Q$ are not edges of $H\setminus X_1$, and so they
	belong to $X_1$. Choose $h\in \{1\LL k\}$ such that $Q$ is a subpath of $P_h$, with ends $u,v$ say, where $u$ is earlier than $v$
	in $P_h$. Consequently $u,v\in \{a_1\LL a_n,b_1\LL b_n\}$. Let $u\in \{a_i,b_i\}$ and $v\in \{a_j,b_j\}$. 
	Since $F$ is $2c$-separated, not both $u=a_i$ and $v=a_j$, and similarly not both  $u=b_i$ and $v=b_j$.
	So either $u=a_i$ and $v=b_j$, or $u=b_i$ and $v=a_j$.
	
	Suppose first that $u=a_i$ and $v=b_j$. Since $Q=P_h[a_i, b_j]$ has length at most $c$, and $R_j=P_h[a_{j+1},b_j]$ has length at least $c+1$, and both 
	$a_i, a_{j+1}$ are earlier than $b_j$, it follows that $R_j$ contains $Q$, and similarly $R_{i-1}$ contains $Q$. In particular, both
	$R_j,R_{i-1}$ contain the end-edges of $Q$, which are in $X_1$, and so $j=i-1$. 
	But then $Q=R_j$ and so has length more than $c$, a contradiction.
	
	Finally, suppose that $u=b_i$ and $v=a_j$. Since $u\ne v$ and $F$ is $2c$-separated, it follows that $v\notin \{b_1\LL b_n\}$,
	and so $a_jb_j$ is the only pair in $F$ incident with $v$. Let $e$ be the edge of $P_h[u,v]$ incident with $v$. Since
	$e\in X_1$, there exists $i'\in \{1\LL n\}$ such that $e\in E(R_{i'})$, and therefore $i'\ne j-1$. 
	Consequently $b_{i'}$ is in $P_h$ and later than $v=a_j$ (and therefore later than $b_i$) in $P_h$, and so $i'>i$.
	Moreover, 
	$a_{i'+1}$ is in $P_h$ 
	and earlier
	than $v$ in $P_h$. Since $F$ is $2c$-separated, it follows that $P_h[a_{i'+1},a_j]$ has length more than $2c$, and so 
	$a_{i'+1}$ is also earlier than $b_i$, and $P_h[a_{i'+1},b_i]$ has length more than $c$ since $Q=P[b_i,a_j]$ has length at most $c$.
	Since $i'\ge i+1$, this contradicts (1) (taking $u,v$ of (1) to be
	$b_i, a_{i'+1}$ respectively). This proves \ref{getpaths}.~\bbox
	
	Finally, here is another lemma we will need:
	
	\begin{thm}\label{shortcut}
		In the same notation, let $a_1b_1\LL a_nb_n$ be a $c$-augmenting sequence, and let $\mac J$ be a partition of $\{1\LL n\}$.
		Then there is
		a $c$-augmenting sequence $a_1'b_1', a_2'b_2'\LL a_m'b_m'$ such that 
		\begin{itemize}
			\item for $1\le i'\le m$, there exist $J\in \mac J$ and $i,j\in J$
			such that $a'_{i'}=a_i$ and $b'_{i'}=b_j$;
			\item for each $J\in \mac J$ there is at most one $i\in J$ such that $a_i\in \{a_1'\LL a_m'\}$, and (therefore) at most one
			$j\in J$ such that $b_j\in \{b_1'\LL b_m'\}$.
		\end{itemize}
	\end{thm}
	\Proof We proceed by induction on $n$. We may assume all members of $\mac J$ are nonempty. If they are all of size one, the result is true,  
	so we may assume that $J_1\in \mac J$ has size at least two.
	Choose $i, j\in J_1$ respectively 
	minimum and maximum; then
	$$a_1b_1\LL a_{i-1}b_{i-1},a_ib_j,a_{j+1}b_{j+1}\LL a_nb_n$$
	is a $c$-augmenting sequence. 
	Let $m=n+i-j$, and define:
	\begin{align*}
		a_h'&=a_h \text{ for } 1\le h\le i\\
		b_h'&=b_h \text{ for } 1\le h\le i-1\\
		a_h'&= a_{h+j-i} \text{ for } i+1\le h\le m\\
		b_h'&= b_{h+j-i} \text{ for } i\le h\le m
	\end{align*}
	Define $f(J_1)=\{i\}$, and for each $J\in \mac J\setminus \{J_1\}$, define
	$$f(J)=\{h:1\le h\le i-1 \text { and } h\in J\}\cup \{h:i+1\le h\le m \text{ and } h+j-i\in J\}.$$
	Then 
	$$\{f(J):J\in \mac J \text{ and } f(J)\ne \emptyset\}$$ 
	is a partition of $\{1\LL m\}$, and the result 
	follows from the inductive hypothesis applied to this partition and $a_1'b_1'\LL a_m'b_m'$. This proves \ref{shortcut}.~\bbox

	\section{The main proof}
	Now we prove \ref{border2}, which we restate for convenience. 
	\begin{thm}\label{treeform2}
		Let $k\ge 0$, $\delta\ge 5$, and $c\ge 1$ be integers, and let $G$ be a graph with no 
		subgraph that is a subdivision of $H_{\delta}^+$. Let $S,T\subseteq V(G)$. Then either 
		\begin{itemize}
			\item there are $k+1$ paths between $S,T$, pairwise at distance greater than $c$; or
			\item there is a set $X\subseteq V(G)$
			with $|X|\le k$ such that every path between $S,T$ contains a vertex at distance at most $(k\delta)^{4k\delta}c$ from some member of $X$.
		\end{itemize}
	\end{thm}
	\Proof We proceed by induction on $k$; the result is trivial for $k=0$, so we assume that $k\ge 1$.
	Define $d = \lceil 2\delta \log_2(k\delta) \rceil$. Then $2^{d-1}/k >(2d)^{\delta-1}$.
	
	
	


	
	
	
	
	

	Let $G$ be a graph with no subgraph that is a subdivision of $H_{\delta}^+$, and let $S,T\subseteq V(G)$.
	We assume 
	\\
	\\
	(1) {\em There is no $X$ with $|X|\le k$,
		such that every path between $S,T$ contains a vertex at distance at most $(k\delta)^{4k\delta}c$ from some member of $X$.}
	\\
	\\
	We must therefore show that there are $k+1$ paths between $S,T$, pairwise at distance more than~$c$. 
	Define $r= d(4d+14) 2^dc$. 
	An $S$-$T$ path $P$ is {\em near-geodesic} if for all $u,v\in V(P)$, either $\dist_P(u,v)\le r$
	or $\dist_G(u,v)>(4d+14)c$. 
	We claim that 
	\\
	\\
	(2) {\em There are $k$ near-geodesic $S$-$T$ paths in $G$, pairwise at distance more than $(4d+14)c$.}
	\\
	\\
	From the inductive hypothesis, since 
	$$(k\delta)^{4k\delta}\ge (4d+14)(2^{d}+1) ((k-1)\delta)^{4(k-1)\delta},$$
	there are $k$ $S$-$T$ paths $Q_1\LL Q_k$,
	pairwise at distance more than $(4d+14)(2^{d}+1)c$.
	
	Let $Z =V(Q_1\cupcup Q_k)$.
	Since $G$ contains no subdivision of $H_\delta^+$, and $2^{d-1}/k >(2d)^{\delta-1}$, \ref{gathering} implies that there is no 
	subgraph of $G$ that is a $Z$-leaved 
	subdivision of $H_{d}$.
	By \ref{newkey}, taking $\ell=(4d+14)c$,
	there exists $Y\supseteq Z$ such that 
	\begin{itemize}
		\item  every vertex in $Y$ has distance at most $2^{d-1}(4d+14)c$ from $Z$;
		\item for all $u,v\in Y$, if $\dist_G(u,v)\le (4d+14)c$, then $\dist_{G[Y]}(u,v)\le d2^{d}(4d+14)c=r$.
	\end{itemize}
	For $1\le i\le k$, there is a path in $G[Y]$ between the ends of $Q_i$, since $Q_i$ is such a path. Let $P_i$ be a 
	shortest such path. If $u,v\in V(P_i)$ with $\dist_{P_i}(u,v)> r$, then 
	$\dist_{G[Y]}(u,v)> r$,
	and so 
	$\dist_G(u,v)>(4d+14)c$, that is, $P_i$ is near-geodesic, for $1\le i\le k$. 
	
	For each $v\in V(P_i)$, since $v\in Y$, it follows that $v$ has distance at most  $2^{d-1}(4d+14)c$ from $Z$, that is, 
	from some $Q_j$, say $Q(v)$. If $u,v\in V(P_i)$ are adjacent, then 
	$$\dist_G(Q(u),Q(v))\le 2^{d}(4d+14)c+1\le (4d+14)(2^{d}+1)c$$
	and so $Q(u)=Q(v)$ since $Q_1\LL Q_k$ pairwise have distance more than $(4d+14)(2^{d}+1)c$. Since $P_i,Q_i$ have the same ends,
	and so $Q(v)=Q_i$ when $v$ is an end of $P_i$, it
	follows that $Q(v)=Q_i$ for all $v\in V(P_i)$, that is, every vertex in $P_i$ has distance at most $2^{d-1}(4d+14)c$ from $Q_i$.
	Consequently, $P_1\LL P_k$ pairwise have distance more than $(4d+14)(2^{d}+1)c-(4d+14)2^{d}c= (4d+14)c$. This proves (2).
	
	
	\bigskip
	
	Fix $S$-$T$ paths $P_1\LL P_k$, each near-geodesic and pairwise at distance more than $(4d+14)c$, and we may choose them 
	minimal with this property; so each has only one vertex in $S$ and one in $T$, its ends. 
	Let $P_h$ have ends $s_h\in S$ and $t_h\in T$, for $1\le h\le k$. Let $\mac P = \{P_1\LL P_k\}$.
	
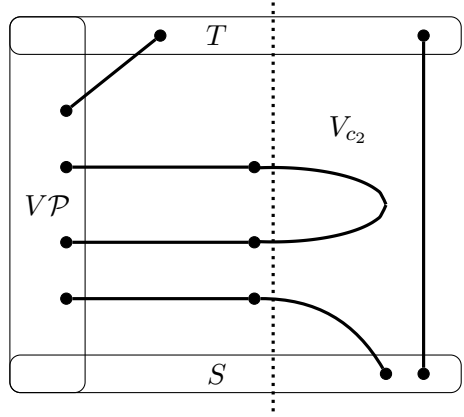
\begin{figure}[h!]
		\centering
		
		\begin{tikzpicture}[scale=.5,auto=left]
			
			\tikzstyle{every node}=[inner sep=1.5pt, fill=black,circle,draw]
			\def\h{5}
			\draw[rounded corners] (-7,-\h) rectangle (-5,\h);
			\draw[rounded corners] (-7,\h) rectangle (5,\h-1);
			\draw[rounded corners] (-7,-\h) rectangle (5,-\h+1);
			\node (s1) at (-.5,1) {};
			\node (s2) at (-.5,-1) {};
			\node (s3) at (-.5,-2.5) {};
			\node (w1) at (-5.5,1) {};
			\node (w2) at (-5.5,-1) {};
			\node (w3) at (-5.5,-2.5) {};
			\node (w4) at (-5.5,2.5) {};
			\node (t4) at (-3,\h-.5) {};
			\node (t3) at (3,-\h+.5) {};
			\node (t5) at (4,-\h+.5) {};
			\node (t6) at (4,\h-.5) {};
			
			\draw[very thick] (w1) to (s1);
			\draw[very thick] (w2) to (s2);
			\draw[very thick] (w3) to (s3);
			\draw[very thick] (w4) to (t4);
			\draw[very thick] (t5) to (t6);
			
			\draw[very thick] (s3) to [bend left=30] (t3);
			
			\def\b{.33}
			\draw[domain=-.5: 3,smooth,variable=\x,very thick] plot ({\x},{\b*sqrt(9-\x*\x)});
			\draw[domain=-.5: 3,smooth,variable=\x,very thick] plot ({\x},{-\b*sqrt(9-\x*\x)});

			\tikzstyle{every node}=[]
			\node at (-1.5,4.5) {$T$};
			\node at (-1.5,-4.5) {$S$};
			\node at (-6,0) {$V\mac P$};
			\node at (2,2) {$V_{(d+3)c}$};
			
			\draw[dotted, very thick] (0,-\h-.5) -- (0,\h+.5);
			
		\end{tikzpicture}
		
		\caption{The four types of leaps. (The thick lines represent paths.)} \label{fig:leaps}
	\end{figure}

	For $p\ge 1$,
	let $V_p$ be the set of vertices at distance more than $p$ from $V\mac P$. Let $L$ be a path of $G$ with ends $a,b$. We say (see Figure \ref{fig:leaps}):
	\begin{itemize}
		\item $L$ is 
		a {\em leap of type 1}
		if $a,b\in V\mac P$, and there exist $x,y\in V(L)$ with $a,x,y,b$ in order, such that  
		the subpaths $L[a,x], L[b,y]$ have length exactly $(d+3)c$, $L[x,y]$ has length at least two,
		and every internal vertex of $L[x, y]$ belongs to $V_{(d+3)c}$. 
		(It follows that $L[a,x]$ is an $(x,V\mac P)$-geodesic, and $L[b,y]$ is a $(y,V\mac P)$-geodesic.)
		\item
		$L$ is a {\em leap of type 2} if $a\in V\mac P$, $b\in (S\cup T)\cap V_{(d+3)c}$, and there exists $x\in V(L)$ such that 
		$L[a,x]$ has length $(d+3)c$, and every  internal vertex of $L[x, b]$ belongs to $V_{(d+3)c}$. 
		\item 
		$L$ is a {\em leap of type 3} 
		if $a\in V\mac P$, 
		$b\in (S\cup T)\setminus V_{(d+3)c}$, and $L$ is a $(b, V\mac P)$-geodesic.
		\item 
		$L$ is a {\em leap of type 4} if $a\in S$ and $b\in T$ and $V(L)\subseteq V_{(d+3)c}$. 
	\end{itemize}
	A {\em leap} is a leap of type 1, 2, 3 or 4.
	
	Let $F$ be the set of all ordered pairs $uv$ such that
	some leap has ends $u,v$. (Thus, if $ab\in F$ then $ba\in F$.)
	\\
	\\
	(3) {\em $F$ is $5r$-jumping in the setting $(S,T,\mac P = \{P_1\LL P_k\})$.}
	\\
	\\
	For $1\le i\le k$ let $Q_i$ be a subpath of $P_i$ of length at most $5r$; thus, $Q_1\LL Q_k$ is a $5r$-barrier
	in the stated setting.
	We may assume (by extending $Q_h$) that for $1\le h\le k$, either $Q_h=P_h$ or $Q_h$ has length exactly $5r$. 
	For $1\le h\le k$, $P_h\setminus V(Q_h)$ has at most two components. If one of them contains $s_h$, call it $A_h$, and otherwise 
	let $A_h$ be the null graph; and if one contains $t_h$ call it $B_h$, and otherwise $B_h$ is null.  
	Choose $q_h\in V(Q_h)$
	for $1\le h\le k$.
	Let $X$ be the set of vertices $v$ of $G$ with $\dist_G(v,A_1\cupcup A_k)\le (d+3)c$ and 
	$\dist_G(v,\{q_1 \LL q_k\})>(k\delta)^{4k\delta}c$; and let $Y$ be 
	the set of $v$ with $\dist_G(v,B_1\cupcup B_k)\le (d+3)c$ and $\dist_G(v,\{q_1 \LL q_k\})>(k\delta)^{4k\delta}c$.
	
	Suppose that $\dist_G(X,Y)\le1$; then $\dist_G(A_1\cupcup A_k, B_1\cupcup B_k)\le 2(d+3)c+1$. Choose $i,j\in \{1\LL k\}$
	such that $\dist_G(A_i, B_j)\le 2(d+3)c+1$. Since 
	$\dist_G(P_i, P_j)>(4d+14)c> 2(d+3)c+1$ for all distinct $i,j$, it follows that $i=j$. Since $\dist_{P_i}(A_i,B_i)\ge 5r+2$ (because both intervals are disjoint from $Q_i$, which has length $5r$), there are vertices $u,v\in P_i$,
	such that $\dist_{P_i}(u,v)\ge 5r+2$ and yet $\dist_G(u,v)\le 2(d+3)c+1$, contradicting that $P_i$ is near-geodesic, 
	since $2(d+3)c+1 \le (4d+14)c$ and $5r+2>r$. This proves that $\dist_G(X, Y)\ge2$, and in particular $X,Y$ are disjoint. 
	
	We claim that for all $v\in V(G)\setminus V_{(d+3)c}$, if $\dist_G(v,\{q_1 \LL q_k\})>(k\delta)^{4k\delta}c$, then 
	every $(v,V\mac P)$-geodesic is either a $(v, A_1\cupcup A_k)$-geodesic or a $(v, B_1\cupcup B_k)$-geodesic.
	Let $J$ be a $(v,V\mac P)$-geodesic,
	and let $b$ be its end in $V\mac P$. Then $b\notin V(Q_1\cupcup Q_k)$ since 
	$\dist_G(v,\{q_1\LL q_k\})>(k\delta)^{4k\delta}c\ge (4d+14)c+ 5r$ and $J$ has length at most $(d+3)c$ and 
	$Q_1\LL Q_k$ all have length at most $5r$. Thus, either $b\in V(A_1\cupcup A_k)$ or $b\in V(B_1\cupcup B_k)$,  and so $J$ is either
	$(v, A_1\cupcup A_k)$-geodesic or a $(v, B_1\cupcup B_k)$-geodesic as claimed. In particular, $v$ belongs to one of $X,Y$.
	
	If $S\cap Y\ne \emptyset$, let $s\in S\cap Y$ and let $J$ be an $(s, V\mac P)$-geodesic, and let $b\in B_1\cupcup B_k$ 
	be the end of $J$ in 
	$V\mac P$.
	Then $J$ is a leap of type 3, and so $sb \in F$ jumps the $5r$-barrier $Q_1\LL Q_k$. Thus,
	we may assume that $S\cap Y=\emptyset$, and similarly $T\cap X=\emptyset$. 
	
	From (1), applied to the set $\{q_1\LL q_k\}$, there is an $S$-$T$ path $P$ in $G$ such that 
	$$\dist_G(P,\{q_1\LL q_k\})>(k\delta)^{4k\delta}c.$$
	Consequently, for each vertex $v\in V(P)\setminus (X\cup Y)$, $\dist_G(v,Q_1\cupcup Q_k)>(k\delta)^{4k\delta}c-5r\ge (d+3)c$, 
	and $\dist_G(v,A_1\cupcup A_k)>(d+3)c$
	(since $v\notin X$), and similarly $\dist_G(v,B_1\cupcup B_k)>(d+3)c$; so $v\in V_{(d+3)c}$, and therefore 
	$V(P)\subseteq X\cup Y\cup V_{(d+3)c}$. Since $P$ has first vertex 
	in $S\cup X$ and last vertex in $T\cup Y$, there is a minimal subpath $Q$ of $P$ between $S\cup X$ and $T\cup Y$, with ends 
	$x\in S\cup X$
	and $y\in T\cup Y$, say, and therefore with no internal vertex in $X\cup Y\cup S\cup T$.

	If $x\in S\setminus X$, then $x\notin Y$ (since $S\cap Y=\emptyset$), and since 
	$\dist_G(x,Q_1\cupcup Q_k)>(k\delta)^{4k\delta}c-5r\ge (d+3)c$, it follows that $x\in V_{(d+3)c}$. So if both  $x\in S\setminus X$
	and $y\in T\setminus Y$, then $Q$ is a leap of type 4 and $xy\in F$ jumps the $5r$-barrier; so from the 
	symmetry we may assume that 
	$x\in X$, and consequently $x\ne y$ since $X\cap (Y\cup T)=\emptyset$. Let $J_x$ be an $(x, V\mac P)$-geodesic, with ends $x$ and $a\in V(A_1\cupcup A_k)$. Thus, $J_x$ has length $(d+3)c$,
	since $x\in X$ and has a neighbour in $V_{(d+3)c}$. If $y\in T\setminus Y$, then $Q\cup J_x$ is a leap of type 2, and $ay\in F$ jumps 
	the $5r$-barrier. Thus, we may assume that $y\in Y$; let $J_y$ be a $(y, V\mac P)$-geodesic, with ends $y,b$ where $b\in V(B_1\cupcup B_k)$.
	Then $Q$ has length at least two since $\dist_G(X,Y)\ge2$; and so $Q\cup J_x\cup J_y$ is a leap of type 1, and $ab\in F$
	jumps the $5r$-barrier. This proves (3).
	
	\bigskip
	From \ref{separateall}, and \ref{organize},
	there is an $r$-augmenting, $2r$-separated sequence $v_1v_2,v_3v_4\LL v_{2n-1}v_{2n}$ in $F$. Let $W=\{v_2,v_3\LL v_{2n-1}\}$.
	Thus, $W\subseteq V\mac P$.
	For $2\le i,j\le 2n-1$, let us say $i,j$ are {\em mated} if $i\ne j$ and $\dist_{U\mac P}(v_i,v_j)\le r$. It follows that if 
	$i,j$ are mated, then one of $i,j$ is odd and the other is even, because $v_1v_2\LL v_{2n-1}v_{2n}$
	is $2r$-separated; and for the same reason for each $v_i\in W$, $i$ is mated with $j$ for at most one $j\in \{2\LL 2n-1\}$. 
	Not all of $v_1\LL v_{2n}$ need be distinct, but if $i\ne j$ and 
	$v_i=v_j$ then $i,j$ are mated, because of the following.
	\\
	\\
	(4) {\em If $i,j\in \{2\LL 2n-1\}$ are distinct and $\dist_G(v_i,v_j)\le (4d+14)c$ then $i,j$ are mated.}
	\\
	\\
	Suppose that $\dist_G(v_i,v_j)\le (4d+14)c$. Consequently $v_i,v_j\in V(P_h)$ for some $h\in \{1\LL k\}$, since 
	$v_i,v_j\in V\mac P$.
	Since $P_h$ is near-geodesic, $\dist_{P_h}(v_i,v_j)\le r$, and so $i,j$ are mated.
	This proves (4).
	
	\bigskip
	
	Some notation: if $P$ is a path of $G$ and $X\subseteq V(G)$, we write $P[X]$ for $P[V(P)\cap X]$.
	For $1\le i\le n$ choose a leap $L_i$ with ends $v_{2i-1},v_{2i}$. If some $L_i$ has type 4 (and hence $i=n=1$), then
	$P_1\LL P_k, L_i$ are $S$-$T$ paths satisfying the theorem; so we may assume that each $L_i$ has type 
	1, 2 or 3. Thus, $L_1,L_n$ have types 2 or 3, and all the others have type 1. 
	
\begin{figure}[h!]
		\centering
		
		\begin{tikzpicture}[xscale = .95, yscale=.5,auto=left]
			
			\tikzstyle{every node}=[inner sep=1.5pt, fill=black,circle,draw]
			\def\h{6}
			\draw[rounded corners] (-7,-\h) rectangle (-5,\h);
			\draw[rounded corners] (-7,\h) rectangle (10,\h-1);
			\draw[rounded corners] (-7,-\h) rectangle (10,-\h+1);
			\node (s1) at (-.5,-1) {};
			\node (s2) at (-.5,-3) {};
			\node (w1) at (-5.5,-1) {};
			\node (w2) at (-5.5,-3) {};
			
			\node (s1a) at (.5,-1) {};
			\node (s2a) at (.5,-3) {};
			\node (s1b) at (3.5,-1) {};
			\node (s2b) at (3.5,-3) {};
			\node (s1c) at (4.5,-1) {};
			\node (s2c) at (4.5,-3) {};
			
			\draw[-] (s1) -- (s1a);
			\draw[-] (s1b) -- (s1c);
			\draw[-] (s2) -- (s2a);
			\draw[-] (s2b) -- (s2c);
			\draw[very thick] (w1) to (s1);
			\draw[very thick] (w2) to (s2);
			\draw[very thick] (s1a) to (s1b);
			\draw[very thick] (s2a) to (s2b);
			
			\node (w3) at (-5.5,1) {};
			\node (s3) at (-.5,1) {};
			\node (s3a) at (.5,1) {};
			\node (s3b) at (3.5,1) {};
			\node (s3c) at (4.5,1) {};
			\node (t3) at (8,\h-.5) {};
			
			\draw[very thick] (w3) to (s3);
			\draw[very thick] (s3a) to (s3b);
			\draw[-] (s3b) to (s3c);
			\draw[-] (s3) to (s3a);
			\draw[very thick] (s3c) to [bend right=40] (t3);

			\node (w4) at (-5.5,2) {};
			\node (s4) at (-.5,4.07) {};
			\node (s4a) at (.5, 4.47) {};
			\node (t4) at (3,\h-.5) {};
			\draw[very thick] (w4) to (s4);
			\draw[-] (s4) to (s4a);
			\draw[very thick] (s4a) to (t4);

			\node (w5) at (-5.5,3) {};
			\node (t5) at (-2,\h-.5) {};

			\draw[very thick] (w5) to (t5);
			
			
			\def\b{.33}
			\draw[domain=-.5: 3,smooth,variable=\x,very thick] plot ({\x+5},{\b*sqrt(9-\x*\x)-2});
			\draw[domain=-.5: 3,smooth,variable=\x,very thick] plot ({\x+5},{-\b*sqrt(9-\x*\x)-2});

			\tikzstyle{every node}=[]
			\node at (1.5,\h-.5) {$T$};
			\node at (1.5,-\h+.5) {$S$};
			\node at (-6,0) {$V\mac P$};
			\node at (-2,\h+.7) {$\le c$};
			\node at (2,\h+1.4) {$> c$};
			\node at (2,\h+.7) {$\le (d+3)c$};
			\node at (6,\h+.7) {$> (d+3)c$};
			\draw[below] (s2) node []           {$v_{2i-1}'$};
			\draw[above] (s1) node []           {$v_{2i}'$};
			\node at (7.5,2) {$R_n$};
			\draw[-{Stealth}, line width = .5] (7,1.8) to [out=235, in=0, looseness=.7] (.5,.7);
			\draw[-{Stealth}, line width = .5] (7.5,2.4) to [out=60, in=270, looseness=.9] (8.4,\h-.5);
			
			\draw[left] (w1) node []           {$v_{2i}$};
			\draw[left] (w2) node []           {$v_{2i-1}$};
			\node at (8.5,-2) {$R_i$};
			\draw[-{Stealth}, line width = .5] (8.5,-1.5) to [out=90, in=5, looseness=.37] (.5,-.5);
			\draw[-{Stealth}, line width = .5] (8.5,-2.5) to [out=270, in=355, looseness=.37] (.5,-3.5);
			\draw[dotted, very thick] (0,-\h) -- (0,\h+1.5);
			\draw[dotted, very thick] (4,-\h) -- (4,\h+1.5);
			
			\node at (-2,-2.5) {$S_{2i-1}$};
			\draw[-{Stealth}, line width = .5] (-3,-2.5) to (-5.5,-2.5);
			\draw[-{Stealth}, line width = .5] (-1,-2.5) to (3.5,-2.5);
			
			\node at (-2.5,-1.5) {$S''_{2i}$};
			\draw[-{Stealth}, line width = .5] (-3.2,-1.5) to (-5.5,-1.5);
			\draw[-{Stealth}, line width = .5] (-1.8,-1.5) to (-.5,-1.5);
			\node at (2,-1.5) {$S'_{2i}$};
			\draw[-{Stealth}, line width = .5] (1.4,-1.5) to (.3,-1.5);
			\draw[-{Stealth}, line width = .5] (2.7,-1.5) to (3.6,-1.5);
			
			\draw[left] (w4) node []           {$v_{2n-1}$};
			\node at (-1.5,2) {$S''_{2n-1}$};
			\draw[-{Stealth}, line width = .5] (-2.5,2.5) to (-4,4);
			\draw[-{Stealth}, line width = .5] (-1.5,2.5) to (-2,3.5);
			\draw[-{Stealth}, line width = .5] (-2.7,1.8) to (-3.5,1);
			
			\node at (2.5,3) {$S'_{2n-1}$};
			\draw[-{Stealth}, line width = .5] (2.5,3.5) to (1.5,4.8);
			\draw[-{Stealth}, line width = .5] (2.5,2.5) to (2,1);
			
			\node at (5,3.7) {$R_n$};
			\draw[-{Stealth}, line width = .5] (4.5,4) to (1.5,4.8);

		\end{tikzpicture}
		\caption{Definitions of $R_i, S_i, S'_i, S''_i$ and $v_i'$.} \label{fig:paths}
	\end{figure}
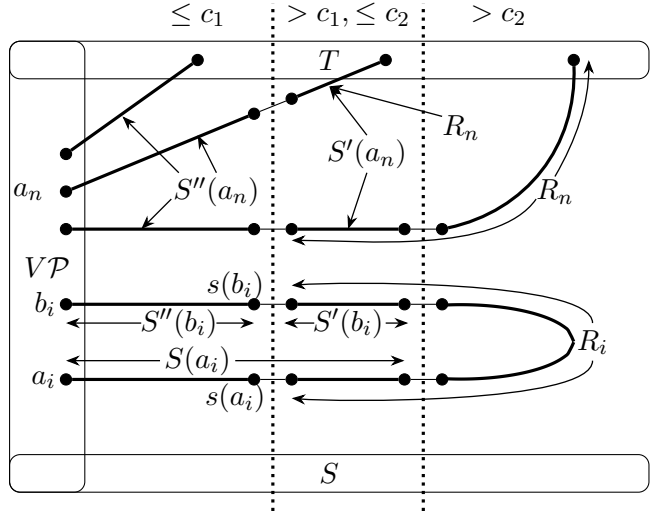

	Let $K$ be the set of all sets $\{i,j\}$ with $i\ne j$ and $2\le i,j\le 2n-1$ such that $i,j$ are mated. 
	For $2\le i\le n$, let $S_{2i-1}$ be the maximal subpath of $L_i$ with one end $v_{2i-1}$ and with length at most $(d+3)c$; 
	and for $1\le i\le n-1$, let $S_{2i}$  be the maximal subpath of $L_i$ with one end $v_{2i}$ and with length at most $(d+3)c$.
	Thus, $S_i$ is defined for $2\le i\le 2n-1$, and 
	$S_i$ has length $(d+3)c$ for $3\le i\le 2n-2$.
	(See Figure \ref{fig:paths}.)
	Let $S'_i=S_i[V_{c}]$ (thus, if $S_i$ has length at most $c$ then $S'_i$ is the 
	null graph).  Let $S''_i=S_i\setminus V_{c}$, and let the ends of $S''_i$ 
	be $v_i, v_i'$.
	For $1\le i\le n$, let 
	$R_i=L_i[V_{c}]$. Thus, $R_i$ is a path unless $L_i$ is a leap of type 3 and has length at most $c$, and then $R_i$ is null. 
	
	We need to be careful with $L_1,L_n$. There are three possibilities for $L_n$ (and the same for $L_1$):
	\begin{itemize}
		\item $L_n$ is a leap of type 2;
		\item $L_n$ is a leap of type 3 and has length more than $c$;
		\item $L_n$ is a leap of type 3 and has length at most $c$.
	\end{itemize}
	Note that, in the second case when $L_n$ has length more than $c$, since $L_n$ is a $(v_{2n},V\mac P)$-geodesic it follows that 
	$V(L_n)\subseteq V(S_{2n-1})\cup V_{c}$, and so 
	$R_n$ joins $v_{2n}$ and a neighbour of $v_{2n-1}'$.

	Let us digress a little, to see the route that lies ahead. We have an $r$-augmenting, $2r$-separated sequence 
	$v_1v_2,v_3v_4\LL v_{2n-1}v_{2n}$, and we could apply \ref{getpaths} to them, and assemble them into $k+1$ paths between $S,T$
	in the graph $H$ of \ref{getpaths}. We could try to convert these to $k+1$ $S$-$T$ paths in $G$ by replacing each pair $v_{2i-1}v_{2i}$
	by the corresponding leap $L_i$ say. Certainly, each path of $H$ is converted into a connected subgraph intersecting both $S,T$. 
	It might not be a path; tracing it gives a walk from $S$ to $T$ that might intersect itself. This is not a problem, we can just shortcut this 
	into a path in the natural way. 
	
	More seriously, while the $k+1$ paths in $H$ are pairwise disjoint, when we convert them to paths of $G$, they 
	might no longer be disjoint, since $L_i, L_j$ might intersect for distinct $i,j$. This is a problem we can 
	easily avoid, by applying \ref{shortcut} to arrange that for each component $D$ of $G[V_{(d+3)c}]$ there is at most one pair $v_{2i-1}v_{2i}$
	for which the corresponding leap passes through $D$. We still need to worry that, say, $S_{2i-1}$ intersects $S_{2j-1}$; but then 
	$2i-1, 2j-1$ would be mated, and hence be joined by a subpath of some $P_h$ of length at most $(4d+14)c$, and therefore would not 
	belong to different paths in the output of \ref{getpaths}, and that problem goes away. 
	
	So we obtain $k+1$ $S$-$T$ paths of $G$, pairwise vertex-disjoint, and it would be nice if they pairwise had distance more than $c$.
	The parts of these paths that are close to $V\mac P$
	are far apart, as we need. The difficulty is that the parts of these paths far away from $V\mac P$ might be too close to each other.
	More exactly, each $L_i$ is composed of sections $S''_{2i-1}, R_i$ and $S''_{2i}$ (let us ignore the first and last leap for simplicity), 
	and the distance between the $S''$ parts of distinct 
	$L_i,L_j$ is satisfactory. Indeed, the parts $S_{2i-1}, S_{2i}$ of $L_i$ are sufficiently far from the parts $S_{2j-1}, S_{2j}$ of 
	$L_j$. But the part of $L_i$ within $V_{(d+3)c}$ (that is, the subpath joining $S_{2i-1}, S_{2i}$) might be uncomfortably close
	to the corresponding part of $L_j$, or to the subpaths $S''_{2j-1}, S''_{2j}$. 
	
	If so, we could choose a short path joining them, and reroute the paths via a process like that of \ref{shortcut}. But we have to repeat this, and while initially these reroutings stay in or close to $V_{(d+3)c}$, when we iterate they can build on one another and become further
	from $V_{(d+3)c}$, and we lose control. The lemma \ref{settle} is designed to handle this.
	
	That ends the digression; let us return to the proof proper.
	For each $\{i,j\}\in K$, if there is a path in $G$ of length at most $c$ between $S'_i,S'_j$, with all vertices in $V_{c}$,
	choose some such path and call it $N_{ij}=N_{ji}$. 
	Let 
	$$\Gamma=\bigcup\left(R_i:1\le i\le  n\right)  \cup \bigcup\left(N_{ij}:\{i,j\}\in K\right).$$
	Thus, $V(\Gamma)\subseteq V_{c}$. 
	\\
	\\
	(5) {\em There is no $V(\Gamma)$-leaved subdivision of $H_{d}$ in $G[V_{c}]$ whose leaves are all in distinct components of $\Gamma$.}
	\\
	\\
	Suppose that $J$ is a $V(\Gamma)$-leaved subdivision of $H_{d}$ in $G[V_{c}]$ whose leaves are all in distinct components of $\Gamma$. Each leaf $u$ of $J$ belongs to a different component $D_u$ say
	of $\Gamma$, and $u$ belongs either to some $V(R_i)\:(1\le i\le n)$ or to some $V(N_{ij})\;(\{i,j\}\in K)$. In either case, 
	there exist $i_u\in \{2\LL 2n-1\}$ with $S_{i_u}'\subseteq D_u$ and a path
	of $D_u$ from $u$ to $S_{i_u}'$; and hence there is a path $Y_u$ from $u$ to $v_{i_u}$, 
	contained in 
	$D_u\cup S_{i_u}$. Let $U$ be the set of leaves of $J$.
	Each path $Y_u$ is from $u$ to $V\mac P$. They might not be vertex-disjoint, but if $Y_{u}$ shares a vertex with $Y_{u'}$,
	this vertex belongs to $S''_{i_{u}}$ and to $S''_{i_{u'}}$ (since $D_u\cap D_{u'}$ is null), and so $\dist_G(v_{i_u},v_{i_{u'}})\le 2c$,
	and hence $i_u,i_{u'}$ are mated by (4). Hence, for each $u\in U$, there is at most one $u'\in U$ different from $u$ such that 
	$Y_u$ intersects $Y_{u'}$; and hence there is a subset $U'\subseteq U$ with $\abs{U'}\ge 2^{d-1}$ such that $Y_u,Y_{u'}$ are vertex-disjoint
	for all distinct $u,u'\in U'$. By \ref{gathering}, since $2^{d-1}/k >(2d)^{\delta-1}$, it follows that $G$ contains a subdivision of $H_{\delta}^+$, a contradiction. This proves (5).

	\bigskip
	
	Choose $t$ maximum such that there is a sequence of paths $M_1\LL M_t$ of $G[V_c]$ satisfying, for each $i\ge 1$:
	\begin{itemize}
		\item $M_i$ has length at most $c$;
		\item the ends of $M_i$ lie in different components  of
		$\Gamma\cup(M_1\cup\cdots\cup M_{i-1})$
		and none of its internal vertices lie in this graph; and
		\item at most one vertex of $M_i$ belongs to $\bigcup_{2\le j\le 2n-1}S'_j$. 
	\end{itemize}
	Let $\Gamma':=\Gamma\cup(M_1\cup\cdots\cup M_t)$. We claim:
	\\
	\\
	(6) {\em Every vertex in $M_1 \cupcup M_t$ has distance at most $(d+1)(c-1)+2$ from $V_{(d+3)c}$.
		Consequently, there is no path $M$ of $G$ (not necessarily of $G[V_c]$) satisfying 
		\begin{itemize}
			\item $M$ has length at most $c$;
			\item the ends of $M$ lie in different components  of
			$\Gamma'$
			and none of its internal vertices lie in this graph; and
			\item at most one vertex of $M$ belongs to $\bigcup_{2\le j\le 2n-1}V(S'_j)$. 
		\end{itemize}
	}
	\noindent
	Let $x\in V(M_1\cupcup M_t)$, and suppose that $\dist_G(x,V_{(d+3)c})>(d+1)(c-1)+2$. If $x$ is not equal or adjacent to any interior vertex of any $M_i$, then $x\in V(\Gamma)$ and there exists $i\in\{1,2,\ldots,t\}$ such that $M_i$ is an edge with one end $x$ and the other end in $V(\Gamma)$. On the other hand,
	if $x$ equals or is adjacent to a vertex $x'$ in the interior of some $M_i$, then $x'\in V(M_1\cupcup M_t)\setminus V(\Gamma)$; and so $\dist_G(x',V_{(d+3)c})>(d+1)(c-1)+1$.
	By (5) and \ref{settle} 
	(taking $\ell=c$ and replacing $G$ by $\Gamma'$ and $F$ by $\Gamma$), either $x'$ is in the interior of some $M_i$ with both ends in $V(\Gamma)$, or there are three components of $\Gamma$
	such that $x'$ has distance at most $(d+1)(c-1)$ from each of them in $\Gamma'$ (and so in $G$).
	
	In summary, one of the following two cases holds:
	\begin{itemize}
		\item there exists $i\in\{1,2,\ldots,t\}$ such that $M_i$ has both ends in $V(\Gamma)$ and contains a vertex of distance more than $(d+1)(c-1)+1$ to $V_{(d+3)c}$ in $G$; or
		
		\item there exists a vertex $x'$ in the interior of some $M_i$, such that $\dist_G(x',V_{(d+3)c})>(d+1)(c-1)$, and $x'$ has distance at most $(d+1)(c-1)$ from three components of $\Gamma$.
	\end{itemize}
	
	In the first case, let $M_i$ have ends $x_1,x_2$. Since $M_i$ has length at most $c$ and contains a vertex of distance more than $(d+1)(c-1)+1\ge c$ to $V_{(d+3)c}$ in $G$,
	it follows 
	that $x_1,x_2\notin V_{(d+3)c}$;
	but then, for $j = 1,2$, $x_j$ belongs either to $S'_{i_j}$
	for some $i_j\in \{2\LL 2n-1\}$, or to $N_{i_jh}$ for some $\{i_j,h\}\in K$. 
	Thus, $\dist_G(v_{i_j},x_j)\le (d+3)c+c$ for $j = 1,2$. Moreover, $S'_{i_1},S'_{i_2}$ belong to different components of $\Gamma_{i-1}:=\Gamma\cup(M_1\cup\cdots\cup M_{i-1})$, and so $i_1\ne i_2$.
	Since
	$M_i$ has length at most $c$, it follows that 
	$$\dist_G(v_{i_1}, v_{i_2})\le 2((d+3)c+c)+c\le (4d+14)c,$$ 
	and so $i_1, i_2$ are mated by (4). 
	We recall that $x_1$ belongs either to $S'_{i_1}$
	or to $N_{i_1h}$ for some $h$. In the second case, $h=i_2$, since $h, i_2$ are both mated with $i_1$. But $S'_{i_1}, S'_{h}$
	belong to the same component of $\Gamma$ (since $N_{i_1h}$ exists), and yet $S'_{i_1}, S'_{i_2}$ belong to different components of 
	$\Gamma_{i-1}$, from the definition of $M_i$, a contradiction. Thus, $x_1$ belongs to $S'_{i_1}$, and similarly $x_2$ belongs to $S'_{i_2}$.
	This contradicts the hypothesis that at most one vertex of $M_i$ is in $\bigcup_{2\le j\le 2n-1}V(S'_j)$.

	Thus, the second case holds; and so there are three components $C_1,C_2,C_3$ of $\Gamma$
	such that $x'$ has distance at most $(d+1)(c-1)$ from each of them. For $j = 1,2,3$, choose $x_j\in C_j$ with distance at most 
	$(d+1)(c-1)$ from $x'$. 
	Since $\dist_G(x',V_{(d+3)c})>(d+1)(c-1)$, it follows that 
	$x_1,x_2,x_3\notin V_{(d+3)c}$.
	Consequently, for $j = 1,2,3$ there exists $i_j\in \{2\LL 2n-1\}$ such that $v_{i_j}$ is at distance at most $(d+3)c+c$ from $x_j$, 
	and  
	$S'_{i_1}, S'_{i_2}, S'_{i_3}$ all belong to different components of $\Gamma$. In particular, $i_1,i_2,i_3$ are all different.
	Therefore, some pair of $i_1,i_2,i_3$ are not mated, say $i_1, i_2$;
	but 
	$$\dist_G(v_{i_1}, v_{i_2})\le 2(d+1)c+2((d+3)c+c)\le (4d+14)c,$$
	contrary to (4). This proves the first statement of (6).
	
	Suppose that there is a path $M$ of $G$ satisfying:
	\begin{itemize}
		\item $M$ has length at most $c$;
		\item the ends of $M$ lie in different components  of
		$\Gamma'$,
		and none of its internal vertices lie in this graph; and
		\item at most one vertex of $M$ belongs to $\bigcup_{2\le j\le 2n-1}V(S'_j)$.
	\end{itemize}
	If some vertex of $M$ lies in $V(M_1\cupcup M_{t})\cup V_{(d+3)c}$, then it has distance at most $(d+1)(c-1)+2$ from $V_{(d+3)c}$;
	and therefore has distance at least $(d+2)c - ((d+1)(c-1)+2)>c$ from $V(G)\setminus V_c$. Consequently $M$ is a path of $G[V_c]$,
	contrary to the maximality of $t$. So $V(M)$ is disjoint from $V(M_1\cupcup M_{t})\cup V_{(d+3)c}$, and hence both its ends lie in
	$V(\Gamma)\setminus V_{(d+3)c}$ and hence in
	$$\bigcup\left(V(S_j'):2\le j\le  2n-1\right)  \cup \bigcup\left(V(N_{ij}):\{i,j\}\in K\right).$$
	At most one end lies in $\bigcup_{2\le j\le 2n-1}V(S'_j)$; so we may assume that one end $x_1$ of $M$ lies in $V(N_{i_1h_1})$
	for some pair $\{i_1,h_1\}\in K$, and its other end $x_2$ belongs either to $V(S_{i_2}')$ for some $i_2\in \{2\LL 2n-1\}$, or to
	$V(N_{i_2h_2})$ for some $\{i_2,h_2\}\in K$. Thus, $\dist_G(x_j, v_{i_j})\le c+(d+3)c$ for $j = 1,2$, and so 
	$$\dist_G(v_{i_1}, v_{i_2})\le 2(c+(d+3)c)+c\le (4d+14)c.$$ 
	By (4) $i_1,i_2$ are equal or mated. They are not equal since $x_1,x_2$ belong to different components of 
	$\Gamma'$,
	so $i_1,i_2$ are mated. By the same argument with $i_1,h_1$ exchanged, $h_1,i_2$ are mated; and so 
	every two 
	of $h_1,i_1,i_2$ are mated, a contradiction. 
	Thus, there is no such $M$. 
	This proves (6).
	
	\bigskip
	
	Let $\mac D$ be the set of components of 
	$\Gamma'$.
	Each $D\in \mac D$ includes at least one of the
	paths $R_j\; (1\le j\le n)$.  For each $D\in \mac D$, let $J_D$ be the set of $j\in \{1\LL n\}$ such that $R_j$ is a non-null subgraph 
	of $D$, and let $I_D$ be the set of all $i\in \{2\LL 2n-1\}$ such that either $i/2\in J_D$ or $(i+1)/2\in J_D$. 
	Let $D^+$ be the union of $D$ and all the paths $S_i$ with $i\in I_D$, and for convenience we define $I_{D^+} = I_D$. 
	(Incidentally, even if $D_1,D_2\in \mac D$ are distinct and therefore disjoint, it is possible that $D_1^+,D_2^+$ might 
	intersect, because there might exist $i_j\in I_{D_j}$ for $j = 1,2$ such that $S''_{i_1},S''_{i_2}$ intersect. 
	But then $i_1, i_2$ would be mated.) 
	Let $\mac D^+$ be 
	the set of all the graphs $D^+$ for $D\in \mac D$, together with either of $L_1, L_{n}$ that has length at most $c$. If $L_1$ has 
	length at most $c$, define $I_{L_1}={2}$, and if $L_n$ has length at most $c$, define $I_{L_n} = \{2n-1\}$. 
	The sets $I_X\;(X\in \mac D^+)$ are nonempty and pairwise disjoint, and their union
	equals $\{2\LL 2n-1\}$.
	
	For each $X\in \mac D^+$, we say $x\in V(X)$ is {\em innocuous in $X$} if 
	there exists $i\in I_X$ such that
	$x\in V(S_i)$, and $S_i[x,v_i]$ has length at most $2c$, and 
	for each vertex  $y$ of $S_i[x,v_i]$, all edges of $D^+$ incident with $y$ belong to $S_i$.
	(Thus, if $L_1$ has length at most $c$, then all vertices of $L_1$ are innocuous in $L_1$, and a similar statement holds if $L_n$ has length at most $c$.)
	We claim:
	\\
	\\
	(7) {\em If $X_1,X_2\in \mac D^+$ are different, and $M$ is a path of length at most $c$ in $G$ between $X_1,X_2$, 
		then for $j = 1,2$, the end of $M$ in $X_j$ is innocuous in $X_j$.}
	\\
	\\
	Suppose first that $\dist_G(M,V\mac P)>c$, and let $M'$ be a minimal 
	subpath of $M$ that has nonempty intersection with two of the graphs in $\mac D^+$ (and therefore with two members of $\mac D$).
	Let the ends of $M'$ be $x_1'\in D_1'$ and $x_2'\in D_2'$, where $D_1', D_2'\in \mac D$.
	From the maximality of $t$ in the definition of $M_1\LL M_t$, it follows that
	at least two vertices of $M'$ belong to $\bigcup_{2\le j\le 2n-1}S'_j$, and since no internal vertex of $M'$ belongs to
	this subgraph, we deduce that $x_1',x_2'$ are both in $\bigcup_{2\le j\le 2n-1}S'_j$.
	Therefore, for $j =1,2$, there exists
	$i_j\in   I_{D_j'}$ such that $x_j'$ belongs to $S'_{i_j}$, with $i_1\ne i_2$. 
	Thus, $\dist_G(v_{i_1}, v_{i_2})\le 2(d+3)c+c\le (4d+14)c$, and so $i_1, i_2$ are mated by (4). Since $i_1\ne i_2$, and 
	$V(M)\subseteq V_c$, it follows that $N_{i_1i_2}$ exists, and so $S'_{i_1},S'_{i_2}$ belong to the same component of $\Gamma'$,
	contradicting that $D_1'\ne D_2'$.  Thus, $\dist_G(M,V\mac P)\le c$.

	Let $M$ have ends $x_j\in X_j$ for $j = 1,2$.
	Since $\dist_G(M,V\mac P)\le c$, it follows that $\dist_G(x_j,V\mac P)\le 2c$, and therefore $\dist_G(x_j, V_{(d+3)c})>(d+1)c$.
	Since $(d+1)c\ge (d+1)(c-1)+2$,
	by (6), $x_j$ is in none of $M_1\LL M_t$. Choose $i_j\in I_{X_j}$ as follows. 
	If $X_j$ is one of $L_1,L_n$  of length at most $c$, let $i_j$ be 2 or $2n-1$ correspondingly.
	If $X_j = D_j^+$ for some $D_j\in \mac D$, then:
	\begin{itemize}
		\item if $x_j\in D_j^+\setminus D_j$, choose $i_j\in I_{D_j}$ with $x_j\in V(S''_{i_j})$;  
		\item otherwise, either there exists $i_j\in I_{D_j}$ with $x_j\in V(S'_{i_j})$, or 
		\item there exist $i_j, h\in I_{D_j}$ such that $i_j, h$ are mated and
		$x_j\in V(N_{i_jh})$.
	\end{itemize}
	Choose $i_j$ as above, for $j = 1,2$. Since $i_j\in I_{X_j}$ for $j = 1,2$, and $ I_{X_1}\cap  I_{X_2}=\emptyset$ since $X_1\ne X_2$,
	we deduce that $i_1\ne i_2$. 
	In each case, it follows that 
	$\dist_G(x_j, v_{i_j})\le (d+3)c+c$, and so
	$\dist_G(v_{i_1},v_{i_2})\le 3c+2(d+3)c\le (4d+14)c$, and so $i_1,i_2$ are mated by (4). Thus, the third bullet above is impossible, as before, and so
	$x_1\in V(S_{i_1})$ and $x_2\in V(S_{i_2})$.

	Since $\dist_G(x_1,V\mac P)\le 2c$, it follows that 
	the subpath of $S_{i_1}$ between $x_1, v_{i_j}$ has length at most $2c$, because it is an $(x_1, V\mac P)$-geodesic. 
	Let $y$ be a vertex of this subpath,
	and let $e$ be an edge of  $X_1$ incident with $y$.
	To show that $x_1$ is innocuous in $X_1$, it remains to show that $e$ is an edge of $S_{i_1}$, for all such $y,e$. 
	Since $\dist_G(y, v_{i_1})\le 2c$, it follows that $\dist_G(y,V_{(d+3)c})>(d+3)c-2c\ge (d+1)(c-1)+2$, and so 
	$y\notin V_{(d+3)c}$ and $y$ belongs to none of $M_1\LL M_t$, by (6). Thus, $e\in E(X_1) \setminus E(M_1\cup\cdots\cup M_t)$, and so either 
	\begin{itemize}
		\item there exists $h_1\in I_{X_1}$ with $e\in E(S_{h_1})$; or
		\item there is a mated pair $h_1,h$ with $h_1,h \in I_{X_1}$ such that $N_{h_1h}$ exists and $e$ is an edge of $N_{h_1h}$.
	\end{itemize}
	In either case $h_1\in I_{X_1}$, and $\dist_G(v_{h_1},v_{i_1})\le 2c+(c+(d+3)c)\le (4d+14)c$. Hence $h_1,i_1$ are equal or mated, by (4). 
	But $i_1,i_2$ are mated, and $h_1\ne i_2$ since $h_1\in I_{X_1}$; and so $h_1=i_1$. 
	Therefore the second bullet above is impossible, and so $e\in E(S_{i_1})$. This proves that
	$x_1$ is innocuous in $X_1$, and similarly $x_2$ is innocuous in $X_2$, and so 
	proves (7).
	
	\bigskip

	We recall that for each $D\in \mac D$, $J_D$ is the set of $j\in \{1\LL n\}$ such that $R_j$ is a non-null subgraph of $D$.
	We would like to apply \ref{shortcut} to the set of sets $\{J_D:D\in \mac D\}$, but it might not be a partition of $\{1\LL n\}$.
	Certainly its union contains $\{2\LL n-1\}$, but we have to be careful about $1,n$.
	There is no $D\in \mac D$ with $1\in J_D$ if and only if $L_1$ is a leap of type 3 of length at most $c$; and the same for $n, L_n$. 
	Let $\mac J$
	be the partition of $\{1\LL n\}$ formed by the sets $\{J_D:D\in \mac D\}$, together with $\{1\}$ if $L_1$ is a leap of type 3 of length at most $c$,
	and $\{n\}$ if $L_n$ is a leap of type 3 of length at most $c$. 
	The sequence $v_1v_2\LL v_{2n-1}v_{2n}$ is $r$-augmenting  and $2r$-separated; and by 
	applying \ref{shortcut} to the partition $\mac J$ and this sequence, we deduce that there is an
	$r$-augmenting, $2r$-separated sequence $w_1w_2\LL w_{2m-1}w_{2m}$ such that, writing $T_j=S_i$ and $T'_j=S'_i$ if $w_j=v_i$, we have:
	\begin{itemize}
		\item $w_1,w_3\LL w_{2m-1}\in \{v_1,v_3\LL v_{2n-1}\}$, and $w_2,w_4\LL w_{2m}\in \{v_2,v_4\LL v_{2n}\}$;
		\item for $1\le i\le m$, either:
		\begin{itemize}
			\item $T'_{2i-1}\cup  T'_{2i}$ is non-null, and there exists $D_i\in \mac D$ such that $V(T'_{2i-1}), V(T'_{2i})\subseteq D_i$; or
			\item $i=1$, and $L_1$ is a leap of type 3 with length at most $c$, and $(w_1,w_2) = (v_1,v_2)$, or
			\item $i=m$, and $L_n$ is a leap of type 3 with length at most $c$,  and $(w_{2m-1}, w_{2m}) = (v_{2n-1},v_{2n})$;
		\end{itemize}
		and 
		\item $D_2\LL D_{m-1}$ and  (if they exist) $D_1,D_m$ are all different.
	\end{itemize}
	To see this, observe that $w_1\in S\setminus V\mac P$, and $w_1\in \{v_1,v_3\LL v_{2n-1}\}$, and therefore $w_1=v_1$, and so if 
	$\{1\}\in \mac J$ then $(w_1,w_2) = (v_1,v_2)$; and similarly if $\{n\}\in \mac J$ then $(w_{2m-1},w_{2m}) = (v_{2n-1},v_{2n})$.
	For $2\le j\le 2m-1$, define $T_j'' = S_i''$ and $w_j'=v_i'$, where $v_i=w_j$.
	Thus, for $2\le i\le 2m-1$, $T''_i$ is the subpath of $T_i$ between $w_i$ and $w_i'$, of length $c$ unless $T_i$ has 
	length less than $c$. Let $w_1' = v_1=w_1\in S$, let $T''_1$ be the one-vertex graph with vertex $w_1'$, let $w_{2m}'=v_{2n}=w_{2m}\in T$,
	and let $T''_{2m}$ be the one-vertex graph with vertex $w_{2m}'$.
	
	If $i, j$ are mated, and $w_{i'}=v_{i}$ and $w_{j'}=v_{j}$, we say $i',j'$ are {\em checkmated}.
	For $1\le i\le m$, if $D_i$ exists (which it does unless $i\in \{1, m\}$), then both $w_{2i-1}', w_{2i}'$ belong to or have a neighbour in $D_i$;
	let $Q_i$ be a path between 
	$w_{2i-1}', w_{2i}'$ with interior in $V(D_i)$. If $D_1$ does not exist, then $L_1$ has length at most $c$ 
	and joins
	$w_1$ and $w_2$; in this case let $Q_1$ be the one-vertex graph with vertex $w_1$. 
	Similarly if $D_m$ does not exist let $Q_m$ be the one-vertex graph with vertex $w_{2m}$.  Thus, for $1\le i\le m$, $T''_{2i-1}\cup Q_i\cup T''_{2i}$ is a path between
	$w_{2i-1},w_{2i}$. 
	\\
	\\
	(8) {\em For all distinct $i,j\in \{1\LL m\}$, if the distance in $G$ between $T''_{2i-1}\cup Q_i\cup T''_{2i}$ and $T''_{2j-1}\cup Q_j\cup T''_{2j}$
		is at most $c$, then one of $2i-1, 2i$ is checkmated with one of $2j-1, 2j$.}
	\\
	\\
	Let $M$ be a path of length at most $c$ with ends $x_1,x_2$, where $x_1\in V(T''_{2i-1}\cup Q_i\cup T''_{2i})$ and 
	$x_2\in V(T''_{2j-1}\cup Q_j\cup T''_{2j})$. Suppose first that $D_i$ exists. By (7), $x_1$ is innocuous in $D_i^+$. 
	Choose $h_1\in I_{D_i}$ with $x_1\in V(T_{h_1})$, with $h_1\in \{2i-1,2i\}$ if possible. Since $T''_{2i-1}\cup Q_i\cup T''_{2i}$ is a path in $D_i^+$ containing
	$x_i$ with both ends in $W\cup \{w_1,w_{2m}\}$, and for each vertex  $y$ of $T_{h_1}[x_1,w_{h_1}]$, all edges of $D_i^+$ incident 
	with $y$ belong to 
	$T_{h_1}$, it follows that $T_{h_1}[x_1,w_{h_1}]$ is a subpath of $T''_{2i-1}\cup Q_i\cup T''_{2i}$, and therefore one of $2i-1, 2i$
	is also a valid choice for $h_1$; and therefore $h_1$ is one of 
	$2i-1,2i$. Moreover, $T_{h_1}[x_1,w_{h_1}]$ has length at most $2c$, and so $\dist_G(x_1, w_{h_1})\le 2c$. 
	
	Now suppose that $D_i$ does not exist. Then 
	$i\in \{1,m\}$; suppose first that $i=1$. Hence $L_1=M_1=T''_2$ has length at most $c$, and contains $x_1$; let $h_1=2$.
	Then $\dist_G(x_1,w_{h_1})\le c$. Similarly, if $i=m$, let $h_1=2m-1$ and it follows that $\dist_G(x_1,w_{h_1})\le c$.
	Thus, whether $D_i$ exists or not, there exists $h_i\in \{2i-1,2i\}$ such that $\dist_G(x_1,w_{h_1})\le 2c$. Define $h_2$
	similarly for $x_2$; then $h_1\ne h_2$ (since $i\ne j$), and it follows that $\dist_G(w_{h_1}, w_{h_2})\le 2(2c) + c\le (4d+14)c$, and so $h_1,h_2$ are checkmated by (4). This proves (8).

	\bigskip
	
	But now the result follows from \ref{getpaths} applied to $w_1w_2\LL w_{2m-1}w_{2m}$, replacing each pair $w_{2i-1}w_{2i}$ in the resulting paths by 
	$T''_{2i-1}\cup Q_i\cup T''_{2i}$.  Let us see this in more detail. Let 
	$$F=\{w_1w_2\LL w_{2m-1}w_{2m}\},$$ 
	and let $H$ be the graph obtained 
	from $U\mac P$ by adding the remainder of $S\cup T$ as vertices, and the ordered pairs in $F$ as (undirected) edges.
	Since 
	$F$ is $r$-jumping (by \ref{organize}) and $2r$-separated, we deduce from \ref{getpaths} that 
	there exist $k+1$ vertex-disjoint $S$-$T$ paths $Z_1\LL Z_{k+1}$ in $H$,
	such that no two of them are joined by a subpath of $U\mac P$ of length at most $r$. Each $Z_s$ is a concatenation
	of subpaths of $U\mac P$ and edges $w_{2i-1}w_{2i}$. 
	
	For $1\le s\le k+1$, let $F_s$ be the set of pairs in $F$ that are edges of $Z_s$. Thus,  $Z_s\setminus F_s$
	is a subgraph of $G$, and each of its components is a subpath of a member of $\mac P$. 
	\\
	\\
	(9) {\em If $a,b\in \{1\LL k+1\}$ are distinct, then $\dist_G(V(Z_a),V(Z_b))> (4d+14)c$.}
	\\
	\\
	Suppose not; then there exist $x\in V(Z_a)$ and $y\in V(Z_b)$ with $\dist_G(x,y)\le (4d+14)c$. 
	Since $x,y\in V\mac P$, at distance at most $(4d+14)c$, 
	both $x,y$ belong to the same member of $\mac P$, say $P_h$.
	Since $\dist_G(x,y)\le (4d+14)c$ and $P_h$ is near-geodesic,
	it follows that $\dist_{P_h}(x,y)\le r$. But $Z_a,Z_b$ are not joined by a subpath of $U\mac P$ of length at most $r$,
	a contradiction.
	This proves (9). 
	
	\bigskip
	
	For each $a\in \{1\LL k+1\}$, let 
	$Y_a$ be the union of $Z_a\setminus F_a$
	and the path $T''_{2i-1}\cup Q_i\cup T''_{2i}$ for each pair $w_{2i-1}w_{2i}\in F_a$. Then $Y_a$ is a connected subgraph of $G$, containing a
	vertex in $S$ and a vertex in $T$.
	\\
	\\
	(10) {\em $Y_1\LL Y_{k+1}$ pairwise have distance more than $c$.}
	\\
	\\
	Suppose that $a,b\in \{1\LL k+1\}$ are distinct, and there exist $x\in V(Y_a)$ and $y\in V(Y_b)$ such that $\dist_G(x,y)\le c$. 
	By (9), it is not the case that 
	$x\in V(Z_a\setminus F_a)$ and $y\in V(Z_b\setminus F_b)$, so we may assume that $y\notin  V(Z_b\setminus F_b)$.
	Choose $w_{2j-1}w_{2j}\in F_b$ such that $y\in V(T''_{2j-1}\cup Q_j\cup T''_{2j})$.
	
	Suppose that $x\notin V(Z_a)$. Then $x\in V(T''_{2i-1}\cup Q_i\cup T''_{2i})$
	for some $w_{2i-1}w_{2i}\in F_a$.
	From (8), one of $2i-1,2i$ (say $i'$) is checkmated with one of $2j-1,2j$ (say $j'$). Hence $w_{i'},w_{j'}$
	belong to the same member
	of $\mac P$, say $P_h$, and $\dist_{P_h}(w_{i'},w_{j'})\le r$. Yet $w_{i'}\in V(Z_a)$ and $w_{j'}\in V(Z_b)$, contradicting that
	$Z_a,Z_b$ are not joined by a subpath of $U\mac P$ of length at most $r$.
	
	So $x\in V(Z_a)$.
	Since $x\in V\mac P$, $\dist_G(x,Q_j)> c$ and so $y\notin V(Q_j)$; and so $y\in V(T''_h)$
	for some $h\in \{2j-1,2j\}$. Since $T''_h$ is a $(y,V\mac P)$-geodesic, and $x,w_h\in V\mac P$ and $w_h\in V(T''_h)$, it follows 
	that 
	$$\dist_G(y,w_h) = \dist_G(y, V\mac P)\le \dist_G(y,x)\le c,$$ 
	and therefore $\dist_G(x,w_h)\le 2c\le (4d+14)c$; but $x\in V(Z_a)$ and $w_h\in V(Z_b)$, contrary to (9).
	This proves (10). 
	
	\bigskip
	
	From (10), this proves \ref{treeform2}.~\bbox
	\section{Concluding remarks}
	
	What about infinite graphs? We assumed that all our graphs were finite at the start of the paper, but augmenting path arguments 
	work fine in infinite graphs (provided we only want some finite number of paths), and the only place in the proof that we used 
	finiteness was in 
	the choice of $M_1\LL M_t$ with $t$ maximum just before step (6) of the main proof. An easy application of Zorn's lemma
	would do instead, so in fact our theorem works for infinite graphs. (And ``path-width'' is better replaced by ``line-width'' for 
	infinite graphs: see~\cite{AS2} for example.)
	
	And for free, we can get a strengthening to graphs with ``bounded coarse line-width''. A {\em $(p,q)$-line-decomposition} of $G$ is a family $(B_t:t\in L)$ of subsets of $V(G)$, 
	where $L$ is a linearly ordered set, such that 
	\begin{itemize}
		\item $\bigcup_{t\in L} G[B_t] = G$;
		\item for all $t_1,t_2,t_3\in T$, if $t_1\le t_2\le t_3$ (where $\le$ is the linear order on $L$) then 
		$B_{t_1}\cap B_{t_3}\subseteq B_{t_2}$; and
		\item for each $t\in L$, $B_t$ is the union of at most $p$ subsets each with diameter in $G$ at most $q$.
	\end{itemize}
	A class of graphs has {\em bounded coarse line-width} if there are $p,q$ such that every graph in the class has a 
	$(p,q)$-line-decomposition (see~\cite{AS3} for a coarse structural characterization of graphs with bounded coarse line-width).
	
	We showed in~\cite{AS1} that for all $p,q$, there exist $\ell,c$ such that every graph that admits a $(p,q)$-line-decomposition also 
	admits an $(\ell,c)$-quasi-isometry to a 
	graph of line-width at most $p$. (See \cite{AS1} for definitions.)
	So we can strengthen our theorem, since its conclusion is invariant under taking quasi-isometries, and obtain that
	the coarse Menger conjecture is true for graphs in any class with bounded coarse line-width:
	
	\begin{thm}
		Let $k, c,p, q\ge 0$ be integers. Then there exists $\ell\ge 0$, such that for every graph $G$ with a $(p,q)$-line-decomposition, and all $S,T\subseteq V(G)$, either:
		\begin{itemize}
			\item there are $k+1$ paths between $S,T$, pairwise at distance more than $c$; or
			\item there is a set $X\subseteq V(G)$
			with $|X|\le k$ such that every path between $S,T$ contains a vertex at distance at most $\ell$ from some member of $X$.
		\end{itemize}
	\end{thm}
	
	\section*{Remarks}
	For the purpose of open access, the authors have applied a CC BY public copyright licence to
	any author accepted manuscript arising from this submission.

\end{document}